\documentclass[12pt]{amsart}

\usepackage[all]{xy}
\usepackage{amscd}
\usepackage{amssymb}
\usepackage{calrsfs}

\newtheorem{theorem}{Theorem}[section]
\newtheorem{lemma}[theorem]{Lemma}
\newtheorem{corollary}[theorem]{Corollary}
\newtheorem{proposition}[theorem]{Proposition}
\newtheorem{conjecture}[theorem]{Conjecture}

\theoremstyle{definition}
\newtheorem{assumption}[theorem]{Assumption}
\newtheorem{remark}[theorem]{Remark}
\newtheorem{definition}[theorem]{Definition}

\newtheorem{example}[theorem]{Example}

\theoremstyle{remark}

\numberwithin{equation}{section}

\def\A{{\mathbb A}}
\def\F{{\mathbb F}}

\def\Q{{\mathbb Q}}
\def\Z{{\mathbb Z}}
\def\C{{\mathbb C}}
\def\O{{\mathcal O}}
\def\P{{\mathbb P}}
\def\R{{\mathbb R}}
\def\L{{\mathcal L}}
\def\X{{\mathfrak X}}

\def\Fr{\text{\rm Fr}}
\def\Gr{\text{\rm Gr}}

\def\Gal{\text{\rm Gal}}
\def\Grass{\text{\rm Gr}}

\def\Im{\text{\rm Im}\,}
\def\Ker{\text{\rm Ker}\,}

\def\id{\text{\rm id}}

\def\PGL{\text{\rm PGL}}

\def\pr{\text{\rm pr}}

\begin{document}

\title[Weight-monodromy conjecture]
{Weight-monodromy conjecture
for $p$-adically uniformized varieties}

\author[Tetsushi Ito]{Tetsushi Ito}
\address{Department of Mathematics, Faculty of Science, Kyoto University,
Kyoto 606-8502, Japan}
\email{tetsushi\char`\@math.kyoto-u.ac.jp}



\subjclass{Primary: 11G25; Secondary: 11F80, 11F85, 14G20, 14F20}
\date{\today}

\begin{abstract}
The aim of this paper is to prove
the weight-monodromy conjecture
(Deligne's conjecture on the purity of
monodromy filtration)
for varieties $p$-adically uniformized
by the Drinfeld upper half spaces of any dimension.
The ingredients of the proof are to prove a special case of
the Hodge standard conjecture,
and apply a positivity argument of Steenbrink, M. Saito
to the weight spectral sequence of Rapoport-Zink.
As an application, by combining our results with
the results of Schneider-Stuhler,
we compute the local zeta functions
of $p$-adically uniformized varieties
in terms of representation theoretic invariants.
We also consider a $p$-adic analogue by using
the weight spectral sequence of Mokrane.
\end{abstract}

\maketitle

\section{Introduction}
\label{SectionIntroduction}

Let $K$ be a complete discrete valuation field
with finite residue field $\F_q$ of characteristic $p>0$,
$\O_K$ the ring of integers of $K$,
and $l$ a prime number different from $p$.
Let $X$ be a proper smooth variety over $K$, and
$V := H^w_{\text{\rm \'et}}(X_{\overline{K}},\Q_l)$
the $l$-adic cohomology
of $X_{\overline{K}} := X \otimes_{K} \overline{K}$ on which
the absolute Galois group $\Gal(\overline{K}/K)$ acts.
We define the inertia group $I_K$ of $K$ by the exact sequence:
$$
\begin{CD}
1 @>>> I_K @>>> \Gal(\overline{K}/K) @>>> \Gal(\overline{\F}_q/\F_q) @>>> 1.
\end{CD}
$$
The monodromy filtration $M_{\bullet}$ on $V$
is determined by the action of $I_K$ on $V$,
and the weight filtration $W_{\bullet}$ on $V$
is determined by the action of a lift of
the geometric Frobenius element
$\Fr_q \in \Gal(\overline{\F}_q/\F_q)$,
which is the inverse of the $q$-th power map
on $\overline{\F}_q$
(for details, see Definition \ref{DefMonodromyFiltration},
 Definition \ref{DefWeightFiltration}).

The weight-monodromy conjecture claims the coincidence
of these two filtrations up to some shift
(\cite{HodgeI}, \cite{WeilII}, \cite{Illusie1},
\cite{Illusie2}, \cite{Illusie3}, \cite{Rapoport-Zink1}, \cite{Rapoport}).
Sometimes, the weight-monodromy conjecture
is also called {\em Deligne's conjecture on the purity
of monodromy filtration} in the literature.

\begin{conjecture}[Weight-monodromy conjecture]
\label{WMC}
$$ M_i V = W_{i+w} V \qquad \text{for all}\ i. $$
\end{conjecture}

Let us briefly recall the theory of
$p$-adic uniformization by Tate, Raynaud, Mumford,
Mustafin, Kurihara, and others
(\cite{Mumford}, \cite{Mustafin}, \cite{Kurihara}).
Let $\widehat{\Omega}^{d}_K$ be
the Drinfeld upper half space of dimension $d \geq 1$,
which is a rigid analytic space obtained by removing
all $K$-rational hyperplanes from $\P^{d}_K$.
We have a natural action of $\PGL_{d+1}(K)$
on $\widehat{\Omega}^{d}_K$.
For a cocompact torsion free discrete subgroup
$\Gamma \subset \PGL_{d+1}(K)$,
the rigid analytic quotient
$\Gamma \backslash \widehat{\Omega}^{d}_K$
has a structure of a projective smooth variety
$X_{\Gamma}$ over $K$
(for details, see \S \ref{SectionP-adic}).
The main theorem of this paper is as follows.

\begin{theorem}[Weight-monodromy conjecture for $X_{\Gamma}$]
\label{MainTheorem}
Let $\Gamma \subset \PGL_{d+1}(K)$ be
a cocompact torsion free discrete subgroup.
Then, the weight-monodromy conjecture (Conjecture \ref{WMC})
holds for $X_{\Gamma}$.
\end{theorem}

Note that Conjecture \ref{WMC} was known to hold
if $X/K$ satisfies at least one of the following conditions.
\begin{enumerate}
\item $X$ has a proper smooth model over $\O_K$
  (\cite{SGA4-III}, \cite{WeilI}, \cite{WeilII}).
\item $X$ is a curve or an abelian variety over $K$
  (\cite{SGA7-I},\ IX).
\item $X$ is a surface
  (\cite{Rapoport-Zink1} for the semistable case,
   \cite{deJong} for the general case).
\item $X$ is a certain threefold with strictly semistable reduction
  (\cite{Ito2}).
\item $K$ is of characteristic $p>0$
  (\cite{WeilII}, \cite{Terasoma}, \cite{Ito1}).
\end{enumerate}
Therefore, Theorem \ref{MainTheorem} is new only
in mixed characteristic
although our proof of Theorem \ref{MainTheorem}
is valid even if $K$ is of characteristic $p>0$
(However, see also Remark \ref{MainTheoremCharP}).
In mixed characteristic and in dimension $\geq 3$,
Conjecture \ref{WMC} is still open up to now.
Therefore, Theorem \ref{MainTheorem} gives us
an interesting class of varieties of any dimension
in mixed characteristic for which Conjecture \ref{WMC} holds.

Our proof of Theorem \ref{MainTheorem}
is based on a careful analysis of the weight spectral sequence
of Rapoport-Zink (\cite{Rapoport-Zink1}).
We prove a special case of the Hodge standard conjecture
and use it to follow a positivity argument of Steenbrink, M. Saito
in \cite{Steenbrink}, \cite{MSaito},
where they proved a Hodge analogue of Conjecture \ref{WMC}
by using polarized Hodge structures
(see also \cite{GuillenNavarroAznar}).
Namely, in our proof of Theorem \ref{MainTheorem},
the Hodge standard conjecture plays a role of
polarized Hodge structures.
Note that in \cite{Ito2}, by basically the same idea
as in this paper, the author proved Conjecture \ref{WMC}
for certain threefolds by using the Hodge index theorem
which is nothing but the Hodge standard conjecture 
for surfaces.

The outline of this paper is as follows.
In \S \ref{SectionHodge}, we fix notation and recall some
basic facts about the Hodge standard conjecture and related
theorems.
Note that, in this paper, we mainly consider varieties
whose $l$-adic cohomology groups are generated by algebraic cycles.
We introduce appropriate notation in \S \ref{SectionNotation}.

In \S \ref{SectionBlowUp},
we recall some properties of blowing-up and
prove the Hodge standard conjecture of blow-ups
under certain assumptions by looking at the asymptotic
behavior of the Lefschetz operator and
the cup product pairings
(Proposition \ref{BlowUpHodgeConjProp}).
Such argument seems well-known
at least in the simplest case of a blow-up
at a point on a surface.
However, since the multiplicative structure of
a blow-up is not so simple in higher dimensions,
we need tedious computations to prove it.

In \S \ref{SectionBlowUpProjectiveSpace},
we consider a projective smooth variety
$B^n$ of dimension $n$ over $\overline{\F}_q$ which
is obtained by successive blowing-ups of $\P^n_{\overline{\F}_q}$
along linear subvarieties.
The main result in \S \ref{SectionBlowUpProjectiveSpace}
is a proof of the Hodge standard conjecture
for such varieties for certain choice of
an ample line bundle (Proposition \ref{MainPropLefschetzHodge}).
Although we can, in principle,
compute everything combinatorially about the variety $B^n$,
the proof of Proposition \ref{MainPropLefschetzHodge}
is far from trivial.
We look at the combinatorial structure of $B^n$
very carefully, and the proof proceeds by induction on $n$.
This proof is inspired by McMullen's proof of an analogue
of the hard Lefschetz theorem for non-rational polytopes
(\cite{McMullen}).

In \S \ref{SectionFiltration}, we recall some basic
facts on the weight spectral sequence of Rapoport-Zink.
By a positivity argument of Steenbrink, M. Saito,
we give a proof of a special case of Conjecture \ref{WMC}
assuming the existence
of $\Q$-structures on $l$-adic cohomology and
the Hodge standard conjecture
(Proposition \ref{WMCHodgeProp}).

In \S \ref{SectionP-adic}, we firstly give a proof of
Theorem \ref{MainTheorem}.
The key point is that there is an explicitly
constructed formal scheme $\widehat{\Omega}^d_{\O_K}$
locally of finite type over $\text{Spf}\,\O_K$,
whose associated rigid analytic space is
$\widehat{\Omega}^{d}_K$.
From the construction, we see that $X_{\Gamma}$ has a proper
semistable model $\X_{\Gamma}$ over $\O_K$.
Moreover, if $\X_{\Gamma}$ is strictly semistable,
we see that all irreducible components
of the special fiber of $\X_{\Gamma}$
are isomorphic to the variety $B^d$
in \S \ref{SectionBlowUpProjectiveSpace}.
Hence we can prove the Hodge standard conjecture for
them by Proposition \ref{MainPropLefschetzHodge}.
The same is true for intersections of them.
By Proposition \ref{WMCHodgeProp},
the proof of Theorem \ref{MainTheorem} is complete.

In the rest of \S \ref{SectionP-adic} and
\S \ref{SectionApplicationTateConjecture},
we also give some applications and complements.
By combining Theorem \ref{MainTheorem}
with the results of Schneider-Stuhler
on the $l$-adic cohomology of $X_{\Gamma}$
(\cite{Schneider-Stuhler}),
we give a proof of Schneider-Stuhler's conjecture
on the filtration $F^{\bullet}$ on
$H^d_{\text{\rm \'et}}(X_{\overline{K}},\Q_l)$
(Theorem \ref{SchneiderStuhlerConjecture}).
As a consequence, we prove the following formula
of the local zeta function $\zeta(s,X_{\Gamma})$
of $X_{\Gamma}$:
\begin{align*}
  \zeta(s,X_{\Gamma}) &:=
    \prod_{k=0}^{2d} \det \big( 1 - q^{-s} \cdot \Fr_q
    \,;\,H^{k}_{\text{\rm \'et}}(X_{\Gamma} \otimes_K \overline{K},\,
    \Q_l)^{I_K} \big)^{(-1)^{k+1}} \\
  &= (1-q^{-s})^{\mu(\Gamma) \cdot (-1)^{d+1}}
        \cdot \prod_{k=0}^{d} \frac{1}{1-q^{k-s}},
\end{align*}
where $\mu(\Gamma)$ is the multiplicity of
the Steinberg representation in the representation
of $\PGL_{d+1}(K)$ induced from the trivial character
on $\Gamma$ (Theorem \ref{LocalZetaFunction}).
In particular, $\zeta(s,X_{\Gamma})$ is independent of $l$.
In \S \ref{SubsectionP-adicAnalogue},
we consider a $p$-adic analogue by using
the weight spectral sequence of Mokrane (\cite{Mokrane}).
After proving a $p$-adic analogue of Conjecture \ref{WMC}
for $X_{\Gamma}$, we also prove that the local zeta function
$\zeta_{\text{\rm $p$-adic}}(s,X_{\Gamma})$
defined via $p$-adic Hodge theory coincides with
$\zeta(s,X_{\Gamma})$ above
(Theorem \ref{MainTheoremp-adic},
 Theorem \ref{LocalZetaFunctionP-adic}).
Finally, in \S \ref{SectionApplicationTateConjecture},
we give an application of these results to the Tate conjecture for
varieties over number fields admitting $p$-adic uniformization.

Needless to say, there is a very important theory
so called the theory of $p$-adic uniformization of
Shimura varieties established by
\v Cerednik, Drinfeld for Shimura curves,
and Rapoport-Zink, Varshavsky in higher dimensions
(\cite{Cherednik}, \cite{Drinfeld},
 \cite{Rapoport-Zink2}, \cite{Varshavsky}).
We expect that Theorem \ref{MainTheorem}
has applications to the zeta functions of
Shimura varieties admitting $p$-adic uniformization
(\cite{Rapoport}, \cite{Rapoport-Zink1}).
We also expect that the results of this paper will establish
a special case of the compatibility between
the global and the local Langlands correspondences
(\cite{Harris}, Problem 1).

\begin{remark}
After the first version of this paper was written,
Ehud de Shalit informed that he obtained
a $p$-adic analogue of Conjecture \ref{WMC} for $X_{\Gamma}$
(Theorem \ref{MainTheoremp-adic})
by a completely different method.
His proof relies on a combinatorial result of
Alon-de Shalit about harmonic cochains on the Bruhat-Tits buildings
(\cite{Alon-deShalit}, \cite{deShalit}).
Wayne Raskind informed that
he and Xavier Xarles studied the cohomology of
varieties with \lq\lq totally degenerate reduction"
(\cite{Raskind-Xarles}).
Because of a special case of the Hodge standard conjecture
proved in this paper,
$p$-adically uniformized varieties are examples
of varieties with \lq\lq totally degenerate reduction"
in the sense of Raskind-Xarles.
\end{remark}

\vspace{0.2in}

\noindent
{\bf Acknowledgments.}
In the preparation of this paper, the author benefited from many people.
Especially, the final part of this work was done
in Arbeitsgemeinschaft Algebraische Geometrie
at Universit\"at zu K\"oln in November and December 2002.
The author is pleased to thank Michael Rapoport
for his continuous advice, encouragements, and discussions.
He also informed that Theorem \ref{MainTheorem}
has an application to the Tate conjecture
for varieties admitting $p$-adic uniformization.
His result is included in \S \ref{SectionApplicationTateConjecture}
with his permission.
The author also would like to thank
Gerd Faltings, Michael Harris, Luc Illusie,
Fumiharu Kato, Kazuya Kato, Minhyong Kim,
Barry Mazur, Yoichi Mieda,  Yukiyoshi Nakkajima,
Richard Taylor, and Teruyoshi Yoshida
for invaluable discussions and comments.
Finally, the author is grateful to his advisor
Professor Takeshi Saito for advice and
support during the graduate course at the University of Tokyo.
A part of this work was done during the author's stay
at Universit\"at zu K\"oln in June 2002,
Universit\'e de Paris-Sud in June 2002,
and Max-Planck-Institut f\"ur Mathematik in Bonn
from November 2002 to October 2003.
The author would like to thank all of them for their support
and cordial hospitality.
The author was supported by the Japan Society for the
Promotion of Science Research Fellowships for Young
Scientists.

\section{Review of the Hodge standard conjecture}
\label{SectionHodge}

\subsection{Notation and assumptions}
\label{SectionNotation}

In this paper, a variety means a separated scheme of finite type
over a field which is geometrically reduced.
We do not assume it is geometrically irreducible.
Let $F$ be an algebraically closed field of any characteristic,
$l$ a prime number different from the characteristic of $F$.
Let $X$ be a projective smooth variety over $F$ of dimension $n$.
Hence $X$ is a disjoint union of irreducible varieties of dimension $n$.
Let $Z^k(X)$ be the group of algebraic cycles on $X$
of codimension $k$, and
$\text{cl}^k \colon Z^k(X) \to H^{2k}_{\text{\rm \'et}}(X,\Q_l(k))$
the cycle map for $l$-adic cohomology.
Let $Z^{k}_{\text{num}}(X) \subset Z^k(X)$ be the subgroup
consisting of algebraic cycles which are
numerically equivalent to zero.
It is known that 
$$ N^k(X) :=  Z^k(X) / Z^{k}_{\text{num}}(X) $$
is a finitely generated free $\Z$-module
(\cite{Kleiman}, Lemma 5-2,
 see also Lemma \ref{QstructureLemma} below).
Here we put the following assumption on $X$.

\begin{assumption}
\label{AssumptionCohomology}
$$ \dim_{\Q_l} H^{k}_{\text{\rm \'et}}(X,\Q_l)
     = \begin{cases}
         \text{\rm rank}_{\Z}\,N^{k/2}(X) & \text{\rm if $k$ is even} \\
         0 & \text{\rm if $k$ is odd}.
       \end{cases} $$
\end{assumption}

\begin{remark}
Assumption \ref{AssumptionCohomology} is independent of $l$
because the $l$-adic Betti numbers are independent of $l$
(This is a consequence of the Weil conjecture, see, for example,
 \cite{Katz-Messing}, Corollary 1).
\end{remark}

\begin{lemma}
\label{QstructureLemma}
Under Assumption \ref{AssumptionCohomology},
$\text{\rm cl}^k$ induces an isomorphism:
$$ N^k(X) \otimes_{\Z} \Q_l \cong H^{2k}_{\text{\rm \'et}}(X,\Q_l(k)). $$
\end{lemma}

\begin{proof}
Fix an integer $k$.
Let $V \subset H^{2k}_{\text{\rm \'et}}(X,\Q_l(k))$ be
a $\Q_l$-vector subspace generated by the image of the cycle map
$\text{\rm cl}^{k}$.
Take elements $x_1,\ldots,x_r \in Z^{k}(X)$ such that
$\{ \text{\rm cl}^{k}(x_1),\ldots,\text{\rm cl}^{k}(x_r) \}$
is a basis of $V$ over $\Q_l$.
Then, it is easy to see that the image of the map
$$ Z^{n-k}(X) \longrightarrow \Z^{\oplus r},\qquad
   y \mapsto (y \cdot x_1,\ldots,y \cdot x_r) $$
is isomorphic to $N^{n-k}(X)$, where $\cdot$ denotes
the intersection product.
Hence we have the following inequalities:
$$ \text{rank}_{\Z}\,N^{n-k}(X) \leq r \leq
   \dim_{\Q_l} H^{2k}_{\text{\rm \'et}}(X,\Q_l)
   = \dim_{\Q_l} H^{2n-2k}_{\text{\rm \'et}}(X,\Q_l). $$
Therefore, Assumption \ref{AssumptionCohomology}
implies that the above inequalities are equalities.
Hence $V = H^{2k}_{\text{\rm \'et}}(X,\Q_l(k))$ for all $k$,
and all cohomology classes are generated by algebraic cycles.
In particular, the numerical equivalence
coincides with the $l$-adic homological equivalence.
We have a surjective map
$N^k(X) \otimes_{\Z} \Q_l \to H^{2k}_{\text{\rm \'et}}(X,\Q_l(k))$,
between $\Q_l$-vector
spaces of the same dimension, hence it is
an isomorphism.
\end{proof}

Namely, $N^k(X)$ defines $l$-independent $\Q$-structures
on the $l$-adic cohomology of $X$. Similarly, by using the cycle map
for crystalline cohomology (\cite{Gillet-Messing}, \cite{Gros}),
it is easy to see that $N^k(X)$ also defines $\Q$-structures on
the crystalline cohomology of $X$.

We define
$$ H^k(X) := 
    \begin{cases}
      N^{k/2}(X) \otimes_{\Z} \R & \text{if $k$ is even} \\
      0 & \text{if $k$ is odd},
    \end{cases} $$
and $H^{\ast}(X) := \bigoplus_{k} H^k(X)$.
We consider $H^{\ast}(X)$ virtually as cohomology with
coefficients in $\R$.
The cup product $\cup$ on $H^{\ast}(X)$
is induced by the intersection product.
Most properties of $l$-adic cohomology
such as pullback, Poincar\'e duality, K\"unneth formula,
etc. are valid for $H^{\ast}(X)$.
If $X$ is irreducible, we have $H^{0}(X) = H^{2n}(X) = \R$ as usual.
To avoid confusion, we always put the subscript \lq\lq \'et''
for \'etale cohomology.

A {\em $\Q$-divisor} (resp. {\em $\R$-divisor})
is a formal sum of divisors with coefficients in
$\Q$ (resp. $\R$).
A {\em $\Q$-divisor} (resp. {\em $\R$-divisor})
$L$ is called {\em ample} if $L$ is of the form
$L = a_1 L_1 + \cdots + a_r L_r$ with
$r \geq 1$, $a_1,\ldots,a_r \in \Q_{>0}$ (resp. $\R_{>0}$),
and $L_1,\ldots,L_r$ ample.
If $L$ is an ample $\Q$-divisor, some positive
integral multiple of $L$ is an ample divisor.
But this does not hold for ample $\R$-divisors.

There is a natural map from the group of
$\R$-divisors on $X$ to $H^2(X)$.
In this paper, we do not usually distinguish
an $\R$-divisor, a formal sum of line bundles with
coefficients in $\R$ (i.e. $\R$-line bundle),
and its class in $H^2(X)$.
We can naturally define pullbacks and restrictions
for $\R$-divisors.
To avoid confusion, for an $\R$-divisor $H$, we sometimes
use the notation $\O(H)$ instead of $H$.
For example, $f^{\ast} \O(H)$ (resp. $\O(H)|_X$)
denotes the pullback (resp. restriction)
as a formal sum of line bundles.

\subsection{Hard Lefschetz conjecture}
\label{SectionHardLefschetz}

As in \S \ref{SectionNotation},
let $X$ be a projective smooth variety over $F$
of dimension $n$ satisfying
Assumption \ref{AssumptionCohomology},
and $L$ an ample $\R$-divisor on $X$.
By taking cup product with $L$, we have an $\R$-linear map
$$
\begin{CD}
L \colon H^k(X) @>>> H^{k+2}(X)
\end{CD}
$$
called the {\em Lefschetz operator}.

\begin{conjecture}[Hard Lefschetz conjecture]
\label{HardLefschetzConjecture}
For all $k$, $L^k$ induces an isomorphism:
$$
\begin{CD}
  L^k \colon H^{n-k}(X) @>{\cong}>> H^{n+k}(X).
\end{CD}
$$
\end{conjecture}

If $L$ is a $\Q$-divisor,
Conjecture \ref{HardLefschetzConjecture}
is proved by Deligne in all characteristics (\cite{WeilII}).
If $F = \C$, Conjecture \ref{HardLefschetzConjecture}
holds by transcendental methods (\cite{WeilKahler}).
The general characteristic 0 case follows
from this by Lefschetz principle.
However, in characteristic $p>0$,
Conjecture \ref{HardLefschetzConjecture} seems open
in general.

Assume that Conjecture \ref{HardLefschetzConjecture}
holds for $(X,L)$.
For an integer $k$ with $0 \leq k \leq n$, we define
the {\em primitive part} $P^k(X)$ by
$$ P^k(X) = \Ker \big( L^{n-k+1} \colon H^k(X)
       \longrightarrow H^{2n-k+2}(X) \big). $$
We define $P^k(X) = 0$ for $k < 0$ or $k > n$.
Then, by Conjecture \ref{HardLefschetzConjecture},
we have the {\em primitive decomposition}:
\begin{equation}
\label{PrimitiveDecomposition}
   H^k(X) = \bigoplus_{i \geq 0} L^i P^{k-2i}(X)
          \cong \bigoplus_{i \geq 0} P^{k-2i}(X)
   \quad \text{for all}\ k.
\end{equation}

Let $X_1,\ldots,X_m$ be the irreducible components of $X$.
Then we have a natural summation map
$$
\begin{array}{rcl}
  \displaystyle
   H^{2n}(X) = \bigoplus_{i=1}^{m} H^{2n}(X_i)
       = \bigoplus_{i=1}^{m} \R
   & \longrightarrow & \R \\
   (a_1,\ldots,a_m) & \mapsto & a_1 + \cdots + a_m,
\end{array}
$$
which we denote by $\sigma$.
For even $k$ with $0 \leq k \leq n$,
we define the pairing $\langle,\rangle_{H^k(X)}$
on $H^k(X)$ by the composite of the following maps:
$$
\begin{array}{rcccl}
 \langle,\rangle_{H^k(X)} \colon H^k(X) \times H^k(X)
   & \longrightarrow
   & H^{2n}(X) & \stackrel{\sigma}{\longrightarrow} & \R \\
 (x, y) & \mapsto & (-1)^{k/2} L^{n-k} x \cup y. & &
\end{array}
$$
Otherwise, we define $\langle,\rangle_{H^k(X)}$ to be
the zero pairing. The pairings $\langle,\rangle_{H^k(X)}$
are nondegenerate by Conjecture \ref{HardLefschetzConjecture}
and Poincar\'e duality.

We denote the restriction of $\langle,\rangle_{H^k(X)}$
to $P^k(X)$ by $\langle,\rangle_{P^k(X)}$.
The decomposition (\ref{PrimitiveDecomposition}) is
an orthogonal decomposition with respect to
$\langle,\rangle_{H^k(X)}$.
Therefore, each $\langle,\rangle_{P^k(X)}$ is nondegenerate.
Moreover, for all $k$ with $0 \leq k \leq n$,
$\langle,\rangle_{H^k(X)}$ is isomorphic to
the alternating sum of $\langle,\rangle_{P^{k-2i}(X)}$ for $i \geq 0$:
\begin{equation}
\label{PrimitiveDecomposition2}
   \langle,\rangle_{H^k(X)}
     = \sum_{i \geq 0} (-1)^{i} \langle,\rangle_{P^{k-2i}(X)}
   \quad \text{on} \quad
   H^k(X) \cong \bigoplus_{i \geq 0} P^{k-2i}(X).
\end{equation}

\begin{remark}
\label{RemarkDimensionPrimitivePart}
The primitive decomposition may depend on the choice of $L$.
However, the dimensions of $P^k(X)$ are independent of $L$
because
$$ \dim_{\R} P^k(X) = \dim_{\R} H^{k}(X) - \dim_{\R} H^{k-2}(X). $$
\end{remark}

\subsection{Hodge standard conjecture}
\label{SubsectionHodgeConj}

As in \S \ref{SectionHardLefschetz},
let $X$ be a projective smooth variety over $F$
of dimension $n$ satisfying
Assumption \ref{AssumptionCohomology},
and $L$ an ample $\R$-divisor on $X$.
Assume further that
the hard Lefschetz conjecture
(Conjecture \ref{HardLefschetzConjecture})
holds for $(X,L)$.

\begin{conjecture}[Hodge standard conjecture]
\label{HodgeConj}
For all $k$, the pairing $\langle,\rangle_{P^{k}(X)}$
is positive definite.
\end{conjecture}

\begin{remark}
Usually, Conjecture \ref{HodgeConj} is stated for all projective
smooth varieties $X$ without Assumption \ref{AssumptionCohomology}
($Hdg(X)$ in \cite{Kleiman}).
It is a part of Grothendieck's standard conjectures,
and known to hold in the following cases.
\begin{enumerate}
\item If $F = \C$, it is a consequence of the existence of
polarized Hodge structures on singular cohomology (\cite{WeilKahler}).
The general characteristic $0$ case follows by Lefschetz principle.
\item If $X$ is of dimension $\leq 2$, it is known
as the Hodge index theorem
(\cite{Kleiman}, 5, \cite{Fulton}, Example 15.2.4).
As a corollary, we can prove
the case $k \leq 1$ by taking hyperplane sections.
\end{enumerate}
\end{remark}

\begin{remark}
\label{RemarkHodgeConj}
Note that Conjecture \ref{HodgeConj} is stated for a pair $(X,L)$.
It seems that Conjecture \ref{HodgeConj} for $(X,L)$
does not automatically imply Conjecture \ref{HodgeConj}
for $(X,L')$ with another choice of $L'$
(see also Corollary \ref{IndependenceCor} below).
\end{remark}

\begin{example}
\label{ExampleProjectiveSpace}
For an ample $\R$-divisor $L$ on $\P^n$,
it is clear that Conjecture \ref{HardLefschetzConjecture}
holds for $(\P^n,L)$.
The only nontrivial primitive cohomology is
$P^0(\P^n) = H^0(\P^n) = \R$
and the pairing $\langle,\rangle_{P^0(\P^n)}$ is
$$ P^0(\P^n) \times P^0(\P^n) = \R \times \R \ni (x,y) \mapsto
       L^n \cdot x y \in \R = H^{2n}(\P^n). $$
Since $L^n$ is positive, this pairing is positive definite.
Hence Conjecture \ref{HodgeConj} holds for $(\P^n, L)$
for all ample $\R$-divisors $L$ on $\P^n$.
\end{example}

\subsection{Signatures}

Let notation be as in \S \ref{SubsectionHodgeConj}.
Firstly, we recall the definition of signatures.
Let $V$ be a finite dimensional $\R$-vector space
with a bilinear symmetric perfect pairing $\langle,\rangle$.
Take a basis $\{ e_1,\ldots,e_r \}$ of $V$
and consider the matrix $M = (m_{ij})$
defined by $m_{ij} = \langle e_i, e_j \rangle$.
By Sylvester's law of inertia,
all eigenvalues of $M$ are in $\R$,
and the number of positive eigenvalues minus
the number of negative eigenvalues is
independent of the choice of $\{ e_1,\ldots,e_r \}$.
We denote this number by
$\text{sign} \big( V, \langle,\rangle \big)$
and call it the {\em signature} of
the pairing $\langle,\rangle$ on $V$.

\begin{proposition}
\label{HodgeConjSignProp}
The Hodge standard conjecture
(Conjecture \ref{HodgeConj}) holds for $(X,L)$
if and only if
$$ \text{\rm sign} \big( H^k(X), \langle,\rangle_{H^k(X)} \big)
     = \sum_{i \geq 0} (-1)^i \dim_{\R} P^{k-2i}(X)
       \quad \text{for all $k$ with $0 \leq k \leq n$}. $$
\end{proposition}

\begin{remark}
By Remark \ref{RemarkDimensionPrimitivePart},
the right hand side of Proposition \ref{HodgeConjSignProp}
does not depend on a choice of $L$.
\end{remark}

\begin{proof}
By definition, Conjecture \ref{HodgeConj}
is equivalent to the following equality:
\begin{equation}
\label{HodgeConjSignature}
   \dim_{\R} P^{k}(X)
     = \text{sign} \big( P^{k}(X), \langle,\rangle_{P^{k}(X)} \big)
   \quad \text{for all $k$ with $0 \leq k \leq n$}.
\end{equation}
By (\ref{PrimitiveDecomposition2}),
we see that, for all $k$ with $0 \leq k \leq n$,
$$ \text{sign} \big( H^k(X), \langle,\rangle_{H^k(X)} \big)
  = \sum_{i \geq 0} (-1)^i
    \,\text{sign} \big( P^{k-2i}(X), \langle,\rangle_{P^{k-2i}(X)} \big). $$
Hence (\ref{HodgeConjSignature})
implies the conditions in Proposition \ref{HodgeConjSignProp}.

Conversely, assume the conditions in
Proposition \ref{HodgeConjSignProp}.
We shall prove (\ref{HodgeConjSignature}) by induction on $k$.
The case $k=0$ is trivial, since $P^0(X) \cong \R$
and $L^n$ is positive.
If (\ref{HodgeConjSignature}) is proved for $k < k_0$,
the condition of Proposition \ref{HodgeConjSignProp}
for $k = k_0$ implies
\begin{align*}
&\hspace*{-0.3in}
     \text{sign} \big( H^{k_0}(X), \langle,\rangle_{H^{k_0}(X)} \big) \\
  &= \sum_{i \geq 0} (-1)^i \dim_{\R} P^{k_0-2i}(X) \\
  &= \dim_{\R} P^{k_0}(X)
  + \sum_{i \geq 1} (-1)^i\,\text{sign} \big( P^{k_0-2i}(X),
        \langle,\rangle_{P^{k_0-2i}(X)} \big).
\end{align*}
By (\ref{PrimitiveDecomposition2}),
we have
  $\dim_{\R} P^{k_0}(X)
       = \text{sign} \big( P^{k_0}(X), \langle,\rangle_{P^{k_0}(X)} \big)$.
\end{proof}

Since the signature is locally constant on the set of
non-degenerate quadratic forms over $\R$,
we have the following implications about
the Hodge standard conjecture (Conjecture \ref{HodgeConj})
for a variety with different choices of
ample $\R$-divisors on it.

\begin{corollary}
\label{IndependenceCor}
\begin{enumerate}
\item Let $L, L'$ be ample $\R$-divisors on $X$.
If the hard Lefschetz conjecture
(Conjecture \ref{HardLefschetzConjecture})
holds for $(X,\,tL + (1-t)L')$
for all $t \in \R$ with $0 \leq t \leq 1$,
then the Hodge standard conjecture
(Conjecture \ref{HodgeConj}) for $(X,L)$, $(X,L')$
are equivalent to each other.
\item In particular, if
Conjecture \ref{HardLefschetzConjecture}
holds for $(X,L)$ for all ample $\R$-divisors $L$ on $X$,
then Conjecture \ref{HodgeConj} for $(X,L)$
for all ample $\R$-divisors $L$ on $X$
are equivalent to each other.
\end{enumerate}
\end{corollary}

\subsection{Products of varieties}
\label{SubsectionProductVarieties}

As in \S \ref{SectionNotation},
let $X$ (resp. $Y$) be a projective smooth variety over $F$
satisfying Assumption \ref{AssumptionCohomology},
and $L_X$ (resp. $L_Y$) an ample $\R$-divisor on
$X$ (resp. $Y$).
Then, by the K\"unneth formula,
we have an isomorphism:
$$ H^k(X \times Y) \cong \bigoplus_{i=0}^{k}
     H^i(X) \otimes_{\R} H^{k-i}(Y). $$
Hence $X \times Y$ also satisfies
Assumption \ref{AssumptionCohomology}.
The condition that $L_X, L_Y$ are ample is equivalent
to the condition that
$\pr_1^{\ast} \O(L_X) + \pr_2^{\ast} \O(L_Y)$ is
ample, where $\pr_1 \colon X \times Y \to X$
and $\pr_2 \colon X \times Y \to Y$
are projections
(\cite{HartshorneAmple}, I, \cite{HartshorneAG}, II).

\begin{proposition}
\label{ConjectureProduct}
\begin{enumerate}
\item If the hard Lefschetz conjecture
  (Conjecture \ref{HardLefschetzConjecture})
  holds for $(X,L_X)$, $(Y,L_Y)$,
  then it also holds for
  $\big( X \times Y,\ \pr_1^{\ast} \O(L_X) + \pr_2^{\ast} \O(L_Y) \big)$.
\item Moreover, if the Hodge standard conjecture
  (Conjecture \ref{HodgeConj})
  holds for $(X,L_X)$, $(Y,L_Y)$,
  then it also holds for
  $\big( X \times Y,\ \pr_1^{\ast} \O(L_X) + \pr_2^{\ast} \O(L_Y) \big)$.
\end{enumerate}
\end{proposition}

\begin{proof}
Here we give a proof by using the
representation theory of the Lie algebra ${\mathfrak sl}(2)$
(for details, \cite{WeilII}, 1.6, Bourbaki Lie, VIII, \S 1).
It is a consequence of Conjecture \ref{HardLefschetzConjecture}
that we can construct ${\mathfrak sl}(2)$-representations on
$H^{\ast}(X)$, $H^{\ast}(Y)$ by using the Lefschetz operator
$L_X,\,L_Y$.
Then, the tensor product $H^{\ast}(X) \otimes_{\R} H^{\ast}(Y)$
has a natural ${\mathfrak sl}(2)$-action which corresponds
to $\pr_1^{\ast} \O(L_X) + \pr_2^{\ast} \O(L_Y)$.
This proves the first assertion.
The second assertion follows from
a basic property of certain pairings defined on
${\mathfrak sl}(2)$-representations over $\R$ (Bourbaki Lie, VIII, \S 1).
\end{proof}

\begin{corollary}
If Conjecture \ref{HardLefschetzConjecture}
(resp. Conjecture \ref{HodgeConj}) holds for
$(X, L_X)$, $(Y, L_Y)$ for all ample $\R$-divisors
$L_X, L_Y$ on $X,Y$,
then it also holds for $(X \times Y,\,L_{X \times Y})$
for all ample $\R$-divisors $L_{X \times Y}$ on $X \times Y$.
\end{corollary}

\begin{proof}
We may assume that $X,Y$ are connected.
Then we have $H^2(X \times Y) \cong H^2(X) \oplus H^2(Y)$,
and all ample $\R$-divisors on $X \times Y$ are necessarily
of the form  $\pr_1^{\ast} \O(L_X) + \pr_2^{\ast} \O(L_Y)$
in $H^2(X \times Y)$.
Hence the assertion follows from Proposition \ref{ConjectureProduct}.
\end{proof}

\section{Cohomology of varieties obtained by blowing-ups}
\label{SectionBlowUp}

In this section, we study the Hodge standard conjecture
(Conjecture \ref{HodgeConj}) for varieties obtained by blowing-ups.
After recalling basic facts about blowing-ups,
we show that, under certain assumptions,
Conjecture \ref{HodgeConj} holds for a blow-up
for certain choice of an ample $\R$-divisor
(for details of the computation,
 see \cite{SGA5}, VII, \cite{HartshorneAG}, II, \cite{Fulton}, 6.7).
We use the same notation as in \S \ref{SectionNotation}.

\subsection{Setup}
\label{SubsectionSetupBlowUp}

Let $X$ be a projective smooth irreducible variety over $F$ of
dimension $n$, and $Y_1,\ldots,Y_r \subset X$
mutually disjoint smooth closed irreducible subvarieties of
codimension $d \geq 2$.
Let $i \colon Y = \coprod_{k=1}^{r} Y_k \hookrightarrow X$
be the closed immersion.
Let $f \colon X' \to X$ be the blow-up of $X$ along $Y$,
$Y' \subset X'$ the strict transform of $Y$, and
$j \colon Y' \hookrightarrow X',\ g \colon Y' \to Y$
natural morphisms.
Assume that $X,Y_1,\ldots,Y_r$ satisfy
Assumption \ref{AssumptionCohomology}
in \S \ref{SectionNotation}.
Then, it is easy to see that $X',Y'$ also satisfy
Assumption \ref{AssumptionCohomology}
in \S \ref{SectionNotation}.
This can be checked from the computations
in \S \ref{CohomologyComputation}
of the cohomology of $X'$ and $Y'$.

We have the following cartesian diagram:
$$
\begin{CD}
Y' @>{j}>> X' \\
@V{g}VV @V{f}VV \\
Y @>{i}>> X
\end{CD}
$$
Since $Y$ is a disjoint union of mutually disjoint
smooth closed irreducible subvarieties of codimension $d$,
$g \colon Y' \to Y$ has a structure of $\P^{d-1}$-bundle.
More precisely, $Y' = \P({\check{N}}_{Y/X})$,
where ${\check{N}}_{Y/X}$ is the conormal bundle
of $Y$ in $X$.
Let
$$ \xi = c_1 \big( \O_{Y'}(1) \big)\ \in\ H^2(Y') $$
be the first Chern class of $\O_{Y'}(1)$ of
the $\P^{d-1}$-bundle $g \colon Y' \to Y$.
The restriction of the line bundle $\O_{X'}(Y')$
to $Y'$ itself is isomorphic to $\O_{Y'}(-1)$.

We define a vector bundle $F$ of rank $d-1$ on $Y'$
by the kernel of
a natural map $g^{\ast}({\check{N}}_{Y/X}) \to \O_{Y'}(1)$.
Then we have the following exact sequence
of vector bundles on $Y'$:
$$
\begin{CD}
0 @>>> F @>>> g^{\ast}({\check{N}}_{Y/X}) @>>> \O_{Y'}(1) @>>> 0.
\end{CD}
$$

\subsection{Cohomology of $X'$ and $Y'$}
\label{CohomologyComputation}

For all $k$, we have the following exact sequence
of cohomology groups (\cite{SGA5}, VII, 8.5):
\begin{equation}
\label{CohomologyBlowUp1}
\begin{CD}
0 @>>> H^{k-2d}(Y) @>{\lambda}>> H^{k-2}(Y') \oplus H^k(X)
  @>{\mu}>> H^k(X') @>>> 0,
\end{CD}
\end{equation}
where $\lambda$ and $\mu$ are defined as follows:
$$
\lambda(y) = \big( g^{\ast}(y) \cup c_{d-1}(\check{F}),\ - i_{\ast}(y) \big),
\qquad \mu(y,x) = j_{\ast}(y) + f^{\ast}(x).
$$
Note that $\check{F}$ is the dual vector bundle of $F$
and $c_{d-1}(\check{F}) \in H^{2d-2}(Y')$
is the $(d-1)$-th Chern class of $\check{F}$.
Define $\lambda' \colon H^{k-2}(Y') \oplus H^k(X) \to H^{k-2d}(Y)$
by $\lambda'(y,x) = g_{\ast}(y)$, then we have
$$ \lambda' \circ \lambda(y)
     = g_{\ast} \big( g^{\ast}(y) \cup c_{d-1}(\check{F}) \big)
     = y \cup g_{\ast} \big( c_{d-1}(\check{F}) \big) = y $$
since
$g_{\ast} \big( c_{d-1}(\check{F}) \big)
   = (-1)^{d-1} g_{\ast} \big( c_{d-1}(F) \big)
   = g_{\ast}(\xi^{d-1}) = 1$ (\cite{SGA5}, VII, 8.4.2).
Therefore, we have the following isomorphism of cohomology groups:
\begin{equation}
\label{CohomologyBlowUp2}
  H^k(X') \stackrel{\mu}{\cong}
  \Ker \big( g_{\ast} \colon H^{k-2}(Y') \to H^{k-2d}(Y) \big) \oplus H^k(X)
\end{equation}

We shall compute the right hand side explicitly.
For an integer $m$ with $0 \leq m \leq d-1$,
the map
$H^{k-2m}(Y) \to H^{k}(Y'),
    \ y \mapsto g^{\ast}(y) \cup \xi^m$
is injective because
$$ g_{\ast} \big((g^{\ast}(y) \cup \xi^m) \cup \xi^{d-1-m} \big)
   = y \cup g_{\ast}(\xi^{d-1}) = y. $$
Let $H^{k-2m}(Y) \xi^m \subset H^{k}(Y')$ be the image of this map.
For simplicity, we sometimes write $y\,\xi^m$ instead of
$g^{\ast}(y) \cup \xi^m$.
Note that $g_{\ast} \colon H^{k}(Y') \to H^{k-2(d-1)}(Y)$
satisfies
$$
g_{\ast}(a\,\xi^m) = a \cup g_{\ast}(\xi^m)
= \begin{cases} a & \text{if $m=d-1$} \\ 0 & \text{otherwise.} \end{cases}
$$
On the other hand, we have an isomorphism
(\cite{SGA5}, VII, 2.2.6):
\begin{align*}
H^{k}(Y') &\cong \bigoplus_{m=0}^{d-1} H^{k-2m}(Y) \xi^m \\
  &= H^{k}(Y) \oplus H^{k-2}(Y) \xi
       \oplus \cdots \oplus H^{k-2(d-1)}(Y) \xi^{d-1}.
\end{align*}
Therefore, by (\ref{CohomologyBlowUp2}), we have the following
explicit description of the cohomology groups of $X'$:
\begin{equation}
\label{CohomologyBlowUp3}
   H^{k}(X') \cong
   H^{k-2}(Y) \oplus H^{k-4}(Y) \xi
   \oplus \cdots \oplus H^{k-2-2(d-2)}(Y) \xi^{d-2}
   \oplus H^{k}(X).
\end{equation}

The following property of the cohomology of $X'$
will be used later.

\begin{proposition}
\label{BlowUpRestrictionProp}
Let $x' \in H^k(X')$ be a cohomology class on $X'$
whose restriction to $Y'$ is zero.
Then $x'$ is of the form $x' = f^{\ast}(x)$ for a unique
$x \in H^{k}(X)$.
\end{proposition}

\begin{proof}
The uniqueness is clear because
$f_{\ast}(x') = f_{\ast} \big( f^{\ast}(x) \big) = x$.
For the existence,
we write $x'$ of the form
$x' = x_0 + x_1\,\xi + \cdots + x_{d-2}\,\xi^{d-2} + x$
via (\ref{CohomologyBlowUp3}).
This means that there is a relation
$$ x' = j_{\ast} \big( x_0 + x_1\,\xi + \cdots + x_{d-2}\,\xi^{d-2} \big)
    + f^{\ast}(x). $$
It is enough to show
$x_0 = x_1 = \cdots = x_{d-2} = 0$.
Since the restriction of $\O_{X'}(Y')$ to $Y'$ is isomorphic
$\O_{Y'}(-1)$, we have $j^{\ast} \big( j_{\ast}(1_{Y'}) \big) = -\xi$.
Therefore,
\begin{align*}
   j^{\ast}(x') &= j^{\ast}j_{\ast} \big( x_0 + x_1\,\xi + \cdots
      + x_{d-2}\,\xi^{d-2} \big) + j^{\ast}f^{\ast}(x) \\
   &= - x_0\,\xi - x_1\,\xi^2 - \cdots
      - x_{d-2}\,\xi^{d-1} + g^{\ast}i^{\ast}(x).
\end{align*}
On the other hand, $j^{\ast}(x') = 0$ by assumption.
Since $H^{\ast}(Y')$ is a free $H^{\ast}(Y)$-module with
basis $1,\xi,\ldots,\xi^{d-1}$,
we have $x_0 = x_1 = \cdots = x_{d-2} = i^{\ast}(x) = 0$.
\end{proof}

\subsection{Multiplicative structure of $H^{\ast}(X')$}
\label{SubsectionMultiplicativeStructure}

The multiplicative structure of $H^{\ast}(X')$ can be computed
as follows. Define $\nabla^{k}(X')$ by
$$ \nabla^{k}(X') = H^{k-2}(Y') \oplus H^k(X) $$
and
$\cup' \colon \nabla^{k}(X') \times \nabla^{m}(X')
    \to \nabla^{k+m}(X')$ by
$$
(y'_1,x_1) \cup' (y'_2,x_2)
= \big( - y'_1 \cup y'_2 \cup \xi + y'_1 \cup g^{\ast}i^{\ast}(x_2)
  + (-1)^{km} y'_2 \cup g^{\ast}i^{\ast}(x_1),\ x_1 \cup x_2 \big).
$$
Then the following diagram is commutative
(\cite{SGA5}, VII, 8.6.3):
$$
\begin{CD}
\nabla^{k}(X') \times \nabla^{m}(X') @>{\mu \times \mu}>>
    H^{k}(X') \times H^{m}(X') \\
@V{\cup'}VV @V{\cup}VV \\
\nabla^{k+m}(X') @>{\mu}>> H^{k+m}(X')
\end{CD}
$$
where $\cup$ denotes the cup product on $H^{\ast}(X')$.
Since the horizontal maps $\mu, \mu \times \mu$ are surjective,
we can compute the multiplicative structure of $H^{\ast}(X')$
in terms of the multiplicative structure of
$\nabla^{\ast}(X') := \bigoplus_{k} \nabla^{k}(X')$.

Recall that the multiplicative structure of $H^{\ast}(Y')$
is computed by Chern classes (\cite{SGA5}, VII, 3).
Namely, we have the following isomorphism of
$H^{\ast}(Y)$-algebras
$$ H^{\ast}(Y') \cong H^{\ast}(Y)[T] / \big( T^{d}
     + c_1(\check{N}_{Y/X})\,T^{d-1}
     + c_2(\check{N}_{Y/X})\,T^{d-2}
     + \cdots + c_d(\check{N}_{Y/X}) \big). $$
which sends $\xi$ to $T$.

There is another relation on $H^{\ast}(X')$ coming from the exact
sequence (\ref{CohomologyBlowUp1}).
Namely, the image of $1_Y \in H^0(Y)$ by $\lambda$,
$$ \lambda(1_Y) = \big( c_{d-1}(\check{F}),\ -i_{\ast}(1_Y) \big)
   \in H^{2d-2}(Y') \oplus H^{2d}(X), $$
maps to 0 in $H^{2d}(X')$.
By definition, $i_{\ast}(1_Y) \in H^{2d}(X)$
is nothing but the cohomology class of $Y$.
We denote it by $[Y]$.
Since the coefficient of $\xi^{d-1}$ in the expression of
$c_{d-1}(\check{F})$ is 1 (\cite{SGA5}, VII, 8.4.2),
we can replace
$\xi^{d-1}$ by $\big( \xi^{d-1} - c_{d-1}(\check{F}),\ [Y] \big)$.

In conclusion,
we can compute the cup product on $H^{\ast}(X')$
in terms of the cohomology groups $H^{\ast}(X), H^{\ast}(Y)$,
the Chern classes of $\check{N}_{Y/X}$,
and the cohomology class $[Y] \in H^{2d}(X)$.

\subsection{Ample divisors on $X'$}

\begin{proposition}
\label{AmpleIntervalProp}
Let notation be as in \S \ref{SubsectionSetupBlowUp}.
Let $L$ be an ample $\R$-divisor on $X$.
There exists $\alpha \in \R_{>0}$
such that, for all $\varepsilon \in \R$,
the $\R$-divisor $L'$ on $X'$ of the form
$$ L' = f^{\ast} \O(L) - \varepsilon Y' $$
is ample if and only if $0 < \varepsilon < \alpha$.
\end{proposition}

\begin{proof}
We use the following basic facts.
Although some results quoted below are only
stated for usual ample divisors in the literature
(for example, see \cite{HartshorneAmple}, \cite{HartshorneAG}),
it is easy to prove them for ample $\R$-divisors.
\begin{enumerate}
\item Take an effective curve $C$ on $X'$
which passes through a point of $Y'$ and
intersects with $Y'$ transversally.
By the Nakai-Moishezon criterion of ampleness
(\cite{HartshorneAmple}, I, \S 5),
we have $\big( f^{\ast} \O(L) - \varepsilon Y' \big) \cdot C > 0$.
Since $Y' \cdot C > 0$,
$\varepsilon$ can not be arbitrary large.
\item If $f^{\ast} \O(L) - \varepsilon Y'$ is ample,
then $\varepsilon > 0$.
To see this, take an effective curve $C$ contained
in a fiber of $g \colon Y' \to Y$.
Since $f^{\ast} \O(L) \cdot C = L \cdot f_{\ast} C = 0$
and $Y' \cdot C < 0$, we have $\varepsilon > 0$
by the Nakai-Moishezon criterion of ampleness.
\item We see that $f^{\ast} \O(L) - \varepsilon Y'$
is ample for sufficiently small $\varepsilon \in \R_{>0}$,
because the restriction of
$\O_{X'}(-Y')$ to $Y'$ is isomorphic to $\O_{Y'}(1)$
(\cite{HartshorneAG}, II, 7.10).
\item The set
$I = \big\{ \varepsilon \in \R \mid f^{\ast} \O(L) - \varepsilon Y'
     \ \text{is ample} \big\}$
is an interval,
because if $D,D'$ are ample $\R$-divisors,
then $t D + (1-t) D'$ are also ample
for all $t \in \R$ with $0 \leq t \leq 1$.
\item Moreover, $I$ is an open interval because,
for a $\R$-divisor $D$ and an ample $\R$-divisor $D'$,
$t D + (1-t) D'$ is ample for sufficiently small
$t \in \R_{>0}$.
\end{enumerate}
By combining above facts, we conclude that
$I = \{ \varepsilon \in \R \mid 0 < \varepsilon < \alpha \}$
for some $\alpha \in \R_{>0}$.
\end{proof}

\begin{remark}
Since all $\R$-divisors $L'$ on $X'$ can be written
as $ L' = f^{\ast} f_{\ast} \O(L') - \varepsilon Y' $
for some $\varepsilon \in \R$,
Proposition \ref{AmpleIntervalProp} determines
the set of all ample $\R$-divisors on $X'$
in some sense.
However, in practice, it seems difficult to
compute the upper bound $\alpha$ in
Proposition \ref{AmpleIntervalProp} explicitly.
\end{remark}

\subsection{Hodge standard conjecture for $X'$}

\begin{proposition}
\label{BlowUpHodgeConjProp}
Let notation be as in \S \ref{SubsectionSetupBlowUp}
(Recall that $X,Y,X',Y'$ satisfy
Assumption \ref{AssumptionCohomology}
in \S \ref{SectionNotation}).
Let $L$ be an ample $\R$-divisor on $X$.
Assume that the hard Lefschetz conjecture
  (Conjecture \ref{HardLefschetzConjecture})
and the Hodge standard conjecture
  (Conjecture \ref{HodgeConj})
hold for $(X,L)$ and  $\big( Y,\O(L)|_Y \big)$,
where $\O(L)|_Y$ denotes the restriction of
$\O(L)$ to $Y$ which is also ample.
Then, there exists $\beta \in \R_{>0}$
such that, for all $\varepsilon \in \R$ with
$0 < \varepsilon < \beta$, the $\R$-divisor $L'$ on $X'$
of the form
$$ L' = f^{\ast} \O(L) - \varepsilon Y' $$
is an ample $\R$-divisor on $X'$ for which
Conjecture \ref{HardLefschetzConjecture}
and Conjecture \ref{HodgeConj} hold.
\end{proposition}

\begin{proof}
First of all, by the isomorphism (\ref{CohomologyBlowUp3}),
we observe that
$H^{\ast}(X')$ is isomorphic to
$H^{\ast}(Y \times \P^{d-2}) \oplus H^{\ast}(X)$
as an {\em $\R$-vector space},
although the multiplicative structures are different.
Nevertheless, if we choose appropriate bases of $H^{\ast}(X')$
and $H^{\ast}(Y \times \P^{d-2}) \oplus H^{\ast}(X)$,
and take the limit $\varepsilon \to 0$,
the matrices representing the action of the Lefschetz operator
and the cup product pairings on $H^{\ast}(X')$
converge to those on
$H^{\ast}(Y \times \P^{d-2}) \oplus H^{\ast}(X)$.

Let us make the above observation precise.
Let $\{ x_{\gamma} \}$ (resp. $\{ y_{\delta} \}$)
be the basis of $H^{\ast}(X)$ (resp. $H^{\ast}(Y)$)
which is independent of $\varepsilon$.
For $\varepsilon \in \R_{>0}$,
we see that
$\big\{ \varepsilon^{-(d-2)/2+k}\,y_{\delta}\,\xi^{k} \big\}$
is a $\varepsilon$-dependent basis of
$H^{\ast}(Y) \xi^k$ for $0 \leq k \leq d-2$.
By (\ref{CohomologyBlowUp3}), we have
the following isomorphism:
$$ H^{\ast}(X') \cong
   H^{\ast}(Y) \oplus H^{\ast}(Y) \xi
   \oplus \cdots \oplus H^{\ast}(Y) \xi^{d-2}
   \oplus H^{\ast}(X). $$
such that the following set
\begin{equation}
\label{BasisBlowUp}
   \big\{ \varepsilon^{-(d-2)/2}\,y_{\delta} \big\} \cup
   \big\{ \varepsilon^{-(d-2)/2+1}\,y_{\delta}\,\xi \big\}
   \cup \cdots \cup 
   \big\{ \varepsilon^{-(d-2)/2+(d-2)}\,y_{\delta}\,\xi^{d-2} \big\}
   \cup \{ x_{\gamma} \}
\end{equation}
is a basis of $H^{\ast}(X')$.
By using the basis (\ref{BasisBlowUp}), we have an isomorphism
\begin{equation}
\label{IsomorphismBlowUp}
   H^{k}(X') \cong
   H^{k-2}(Y \times \P^{d-2}) \oplus H^{k}(X)
\end{equation}
as an {\em $\R$-vector space}
which sends $\varepsilon^{-(d-2)/2+m}\,\xi^{m}$
to $c_1 \big( \O_{\P^{d-2}}(1) \big)^{\!m}$.
The reason we multiply $\xi^k$ by $\varepsilon^{-(d-2)/2 + k}$
in (\ref{BasisBlowUp}) is as follows.
If we take the cup product of two elements in (\ref{BasisBlowUp}),
the coefficient of the top degree term $\xi^{d-2}$ is
a constant function in $\varepsilon$.
If $\xi^{k}$ for $k > d-2$ appears in the cup product,
we can replace
$\xi^{d-1}$ by $\big( \xi^{d-1} - c_{d-1}(\check{F}),\ [Y] \big)$
in $H^{\ast}(X')$
(see \S \ref{SubsectionMultiplicativeStructure}).
However, since the order of $\varepsilon$ in
the coefficient of $\xi^{d-1}$ is one,
this term vanishes if we take the limit $\varepsilon \to 0$.
Therefore, by using the basis (\ref{BasisBlowUp}),
we can consider the asymptotic behavior
of the Lefschetz operator and the cup product pairings
on $H^{\ast}(X')$ as $\varepsilon \to 0$.

We put
$\widetilde{H}^k(X') := H^{k-2}(Y \times \P^{d-2}) \oplus H^{k}(X)$
and define the map
\begin{equation}
\label{DirectSumLefschetzOperator}
\begin{CD}
  \widetilde{L} \colon \widetilde{H}^{k}(X')
  @>>> \widetilde{H}^{k+2}(X')
\end{CD}
\end{equation}
as follows. Let $\pr_1 \colon Y \times \P^{d-2} \to Y$,
$\pr_2 \colon Y \times \P^{d-2} \to \P^{d-2}$ be projections.
Then
$$ L_{Y \times \P^{d-2}} = \pr_1^{\ast} \big( \O(L)|_Y \big)
       + \pr_2^{\ast} \O_{\P^{d-2}}(1) $$
is an ample $\R$-divisor on $Y \times \P^{d-2}$.
We define the map $\widetilde{L}$
as a direct sum of the Lefschetz operators
$L_{Y \times \P^{d-2}}$ and $L$.
By Example \ref{ExampleProjectiveSpace},
Proposition \ref{ConjectureProduct},
the hard Lefschetz conjecture
(Conjecture \ref{HardLefschetzConjecture})
and the Hodge standard conjecture
(Conjecture \ref{HodgeConj})
hold for $(Y \times \P^{d-2},\,L_{Y \times \P^{d-2}})$.
Therefore, we see that
$\widetilde{L}^k$ induces the following isomorphism:
$$
\begin{CD}
  \widetilde{L}^k \colon \widetilde{H}^{n-k}(X')
  @>{\cong}>> \widetilde{H}^{n+k}(X').
\end{CD}
$$
Namely, $\widetilde{L}$ satisfies an analogue of
the hard Lefschetz conjecture
(Conjecture \ref{HardLefschetzConjecture}).
Moreover, an analogue of
the Hodge standard conjecture (Conjecture \ref{HodgeConj})
holds for $\widetilde{L}$ in the following sense.
Let $\widetilde{P}^k(X')$ be the kernel of $\widetilde{L}^{n-k+1}$
on $\widetilde{H}^{k}(X')$,
which is isomorphic to $P^{k-2}(Y \times \P^{d-2}) \oplus P^{k}(X)$.
Then the following pairing on $\widetilde{P}^k(X')$
$$ \widetilde{P}^k(X') \times \widetilde{P}^k(X')
     \longrightarrow \widetilde{H}^{2n}(X') = \R,
   \qquad (a,b) \mapsto
     (-1)^{k/2}\,\widetilde{L}^{n-k}\,a\,\widetilde{\cup}\,b $$
is positive definite, where $\widetilde{\cup}$
is a difference of the cup products on
$H^{\ast}(Y \times \P^{d-2})$ and $H^{\ast}(X)$
defined as follows:
$$
\widetilde{\cup} \colon
 \widetilde{H}^{\ast}(X') \times \widetilde{H}^{\ast}(X')
 \ni \big( (y_0,x_0),\ (y_1,x_1) \big) \mapsto - y_0 \cup y_1 + x_0 \cup x_1.
$$
Note that the minus sign is inevitable for
the Hodge standard conjecture (Conjecture \ref{HodgeConj})
because, for $(y,x) \in \widetilde{P}^k(X')$,
the degree of $y$ is smaller than the degree of $x$ by 2.

Therefore, 
to prove Proposition \ref{BlowUpHodgeConjProp},
it is enough to prove the following lemma.

\begin{lemma}
\label{HodgeLemma}
If we take the limit $\varepsilon \to 0$,
the matrices representing
the action of $L' = f^{\ast} \O(L) - \varepsilon Y'$
with respect to the basis (\ref{BasisBlowUp})
converge to the matrices representing
$\widetilde{L}$ via (\ref{IsomorphismBlowUp}).
Similarly, for all $k$, the matrices representing
the cup product pairing between
$H^{k}(X')$ and $H^{2n-k}(X')$
converge to the matrices representing
the pairing $\widetilde{\cup}$ between
$\widetilde{H}^{k}(X')$ and
$\widetilde{H}^{2n-k}(X')$.
\end{lemma}

We shall compute the action of
$L' = f^{\ast} \O(L) - \varepsilon Y'$.
Since $[Y'] = j_{\ast}(1_{Y'}) \in H^2(X')$,
we compute the cup product
\begin{equation}
\label{CupProductLefschetzOperator}
   \big( f^{\ast} \O(L) - \varepsilon j_{\ast}(1_{Y'}) \big)
      \cup \big( f^{\ast}(x) + j_{\ast}(y) \big)
\end{equation}
for $x \in H^k(X),\ y \in H^{k-2}(Y')$ with
$g_{\ast}(y) = 0$
(see (\ref{CohomologyBlowUp2})).
Since (\ref{CupProductLefschetzOperator}) is equal to
$$ f^{\ast} \O(L) \cup f^{\ast}(x)
         + f^{\ast} \O(L) \cup j_{\ast}(y)
         - \varepsilon j_{\ast}(1_{Y'}) \cup f^{\ast}(x)
         - \varepsilon j_{\ast}(1_{Y'}) \cup j_{\ast}(y), $$
we consider these four terms separately.
Firstly, since the pullback preserves the cup product,
we have
$ f^{\ast} \O(L) \cup f^{\ast}(x)
     = f^{\ast} \big( \O(L) \cup x \big)
     = f^{\ast}(L \cup x). $
Secondly, we have
$ f^{\ast} \O(L) \cup j_{\ast}(y)
     = j_{\ast} \big(j^{\ast} f^{\ast} \O(L) \cup y \big)
     = j_{\ast} \big(g^{\ast} i^{\ast} \O(L) \cup y \big)
     = j_{\ast} \big(g^{\ast} (\O(L)|_Y) \cup y \big) $
because $f \circ j = i \circ g$.
Similarly, we have
$ - \varepsilon j_{\ast}(1_{Y'}) \cup f^{\ast}(x)
     = - \varepsilon j_{\ast} \big( 1_{Y'} \cup j^{\ast} f^{\ast}(x) \big)
     = - \varepsilon j_{\ast} \big( g^{\ast} (x|_Y) \big). $
Finally, by the self-intersection formula
$j^{\ast} j_{\ast} (1_{Y'}) = - \xi$,
we have
$ - \varepsilon j_{\ast}(1_{Y'}) \cup j_{\ast}(y)
     = - \varepsilon j_{\ast} \big( j^{\ast} j_{\ast} (1_{Y'}) \cup y \big)
     = j_{\ast}(\varepsilon \xi \cup y). $
Therefore, we conclude that
(\ref{CupProductLefschetzOperator}) is equal to
$$ f^{\ast}(L \cup x)
     + j_{\ast} \Big( \big( g^{\ast} (\O(L)|_Y)
           + \varepsilon \xi \big) \cup y \Big)
     - \varepsilon j_{\ast} \big( g^{\ast} (x|_Y) \big). $$
It is clear that the first term
$f^{\ast}(L \cup x)$ corresponds to the action of $L$
on $H^{\ast}(X)$ via (\ref{IsomorphismBlowUp}).
The second term
$j_{\ast} \big( \big( g^{\ast} (\O(L)|_Y)
     + \varepsilon \xi \big) \cup y \big)$
almost corresponds to the action of $L_{Y \times \P^{d-2}}$
on $H^{\ast}(Y \times \P^{d-2})$ via (\ref{IsomorphismBlowUp}).
The only difference is that $\xi^{d-1}$ may appear
in the cup product
$\big( g^{\ast} (\O(L)|_Y) + \varepsilon \xi \big) \cup y$.
Since the order of $\varepsilon$ in the coefficient of $\xi^{d-1}$
is higher than others, if we replace $\xi^{d-1}$ by
$\big( \xi^{d-1} - c_{d-1}(\check{F}),\ [Y] \big)$,
this contribution vanishes if we take the limit $\varepsilon \to 0$.
The contribution from the third term
$- \varepsilon j_{\ast}\big( g^{\ast} (x|_Y) \big)$
vanishes if we take the limit $\varepsilon \to 0$.
Hence the first assertion of Lemma \ref{HodgeLemma}
is proved.

Finally, we shall compute the cup product pairing on $H^{\ast}(X')$.
According to (\ref{CohomologyBlowUp2}),
we compute the cup product
\begin{equation}
\label{CupProductPairing}
   \big( f^{\ast}(x_0) + j_{\ast}(y_0) \big)
       \cup \big( f^{\ast}(x_1) + j_{\ast}(y_1) \big)
\end{equation}
for
$x_0 \in H^k(X)$,
$y_0 \in H^{k-2}(Y')$,
$x_1 \in H^{2n-k}(X)$,
$y_1 \in H^{2n-k-2}(Y')$
with $g_{\ast}(y_0) = 0,\ g_{\ast}(y_1) = 0$.
Since $(\ref{CupProductPairing})$ is in $H^{2n}(X') = \R$,
and $f_{\ast} \colon H^{2n}(X') \to H^{2n}(X)$ is an isomorphism,
it is enough to compute
\begin{align*}
 f_{\ast} \Big( \big( f^{\ast}(x_0) + j_{\ast}(y_0) \big)
       \cup \big( f^{\ast}(x_1) + j_{\ast}(y_1) \big) \Big)
   &= f_{\ast} \big( f^{\ast}(x_0) \cup f^{\ast}(x_1) \big) \\
   &\hspace*{-2.3in}
    + f_{\ast} \big( f^{\ast}(x_0) \cup j_{\ast}(y_1) \big)
    + f_{\ast} \big( j_{\ast}(y_0) \cup f^{\ast}(x_1) \big)
    + f_{\ast} \big( j_{\ast}(y_0) \cup j_{\ast}(y_1) \big).
\end{align*}
For the first term,
we have
$ f_{\ast} \big( f^{\ast}(x_0) \cup f^{\ast}(x_1) \big) = 
   f_{\ast} \big( f^{\ast}(x_0 \cup x_1) \big) = x_0 \cup x_1 $
since the pullback preserves the cup product and
$f_{\ast} \circ f^{\ast} = \id$.
For the second term, we have
$ f_{\ast} \big( f^{\ast}(x_0) \cup j_{\ast}(y_1) \big)
   = x_0 \cup f_{\ast} \big( j_{\ast}(y_1) \big)
   = x_0 \cup i_{\ast} \big( g_{\ast}(y_1) \big)
   = 0 $
because $f \circ j = i \circ g$ and $g_{\ast}(y_1) = 0$.
Similarly, the third term
$f_{\ast} \big( j_{\ast}(y_0) \cup f^{\ast}(x_1) \big)$ is zero
because $g_{\ast}(y_0) = 0$.
For the last term, by the self-intersection formula
$j^{\ast} j_{\ast} (y_0) = - \xi \cup y_0$,
we have
$ f_{\ast} \big( j_{\ast}(y_0) \cup j_{\ast}(y_1) \big)
    = f_{\ast} \Big( j_{\ast} \big( j^{\ast}j_{\ast}(y_0) \cup y_1 \big) \Big)
    = i_{\ast} \big( g_{\ast}(- \xi \cup y_0 \cup y_1) \big). $
Since $g_{\ast}(\xi^{d-1})=1$ and $g_{\ast}(\xi^{m})=0$
for $0 \leq m \leq d-2$,
if we take the limit $\varepsilon \to 0$,
the matrix representing the pairing
$$ H^{k-2}(Y') \times H^{2n-k-2}(Y') \ni (y_0,y_1)
   \mapsto i_{\ast} \big( g_{\ast}(- \xi \cup y_0 \cup y_1) \big)
   \in H^{2n}(X) = \R $$
with respect to the basis (\ref{BasisBlowUp})
converges to the {\em minus} of the matrix representing
the cup product pairing
$$ H^{k-2}(Y \times \P^{d-2}) \times H^{2n-k-2}(Y \times \P^{d-2})
   \longrightarrow H^{2n-4}(Y \times \P^{d-2}) = \R $$
via (\ref{IsomorphismBlowUp}).
Therefore, 
the second assertion of Lemma \ref{HodgeLemma} is proved.
Now the proof of Proposition \ref{BlowUpHodgeConjProp}
is complete.
\end{proof}

\begin{corollary}
Let notation and assumptions be as in
Proposition \ref{BlowUpHodgeConjProp}.
\begin{enumerate}
\item If the hard Lefschetz conjecture
  (Conjecture \ref{HardLefschetzConjecture})
  holds for $(X',L')$ for all ample $\R$-divisors $L'$ on $X'$
  of the form
  $L' = f^{\ast} \O(L) - \varepsilon Y'\ (\varepsilon \in \R_{>0})$,
  then the Hodge standard conjecture (Conjecture \ref{HodgeConj})
  holds for $(X',L')$ for all ample $\R$-divisors $L'$ on $X'$
  of the form $L' = f^{\ast} \O(L) - \varepsilon Y'$.
\item In particular, if the hard Lefschetz conjecture
  (Conjecture \ref{HardLefschetzConjecture})
  holds for $(X',L')$ for all ample $\R$-divisors $L'$ on $X'$,
  then the Hodge standard conjecture (Conjecture \ref{HodgeConj})
  holds for $(X',L')$ for all ample $\R$-divisors $L'$ on $X'$.
\end{enumerate}
\end{corollary}

\begin{proof}
By Corollary \ref{IndependenceCor},
the assertion easily follows from
Proposition \ref{BlowUpHodgeConjProp}.
\end{proof}

\begin{remark}
Our proof of Proposition \ref{BlowUpHodgeConjProp}
is based on a consideration of the asymptotic behavior
of the cup product pairing as $\varepsilon \to 0$.
Therefore, by our method, we can not prove whether
we can take $\beta = \alpha$.
In other words, we can not prove
the Hodge standard conjecture (Conjecture \ref{HodgeConj})
for {\em all} ample $\R$-divisors $L'$ on $X'$
of the form $L' = f^{\ast} \O(L) - \varepsilon Y'$.
\end{remark}

\section{Some blowing-ups of projective spaces}
\label{SectionBlowUpProjectiveSpace}

In this section, we consider a projective smooth variety
$B^n$ of dimension $n$ over $\overline{\F}_q$ which
is obtained by successive blowing-ups of $\P^n_{\overline{\F}_q}$
along linear subvarieties.
All varieties appearing in this section are defined over
$\overline{\F}_q$ and satisfy Assumption \ref{AssumptionCohomology}
in \S \ref{SectionNotation}.
Hence, for simplicity, we omit the subscript
\lq\lq $\overline{\F}_q$'' and
denote $\P^n_{\overline{\F}_q}$ by $\P^n$, etc.
This does not cause any confusion.
We use the same notation as in \S \ref{SectionNotation}.

\subsection{Construction of $B^n$}
\label{SubsectionConstruction}

Let $n \geq 1$ be an integer. We have a natural map
$\varphi \colon \A^{n+1} \backslash \{ 0 \} \to \P^n$.
A subvariety $V \subset \P^n$ is called
a {\em linear subvariety}
if $\varphi^{-1}(V) \cup \{ 0 \} \subset \A^{n+1}$
is a linear subspace.
Let $\Grass_d(\P^n)$ be the Grassmann variety of
linear subvarieties of dimension $d$ in $\P^n$.
Let $\Grass_{\ast}(\P^n)$ be the disjoint union
of $\Grass_d(\P^n)$ for all $d$.
$\Grass_d(\P^n)(\F_q)$ is the set of
linear subvarieties of dimension $d$ in $\P^n$
defined over $\F_q$.

We construct a birational map $f \colon B^n \to \P^n$ as follows.
Firstly, we put $Y_0 := \P^n$.
Then, let $Y_1$ be the blow-up of $Y_0$ along
the disjoint union of all $\F_q$-rational points on $\P^n$.
Let $Y_2$ be the blow-up of $Y_1$ along
disjoint union of all strict transforms of lines in $\P^n$.
Similarly, we construct projective smooth varieties
$Y_0,\ldots,Y_{n-1}$ inductively as follows.
Assume that $Y_k$ was already constructed.
Let $Z_k \subset Y_k$ be the union of
all strict transforms of linear subvarieties
of dimension $k$ in $\P^n$ defined over $\F_q$.
These are disjoint because all intersections of them
were already blown-up.
Let $g_k \colon Y_{k+1} \to Y_k$ be the blow-up of $Y_k$
along $Z_k$.
Finally, we put $B^n := Y_{n-1}$
and $f := g_{0} \circ \cdots \circ g_{n-2} \colon B^n \to \P^n$.
We have the following sequence of blowing-ups.
$$
\begin{CD}
  B^n = Y_{n-1} @>{g_{n-2}}>> Y_{n-2} @>{g_{n-3}}>> \cdots
    @>{g_1}>> Y_1 @>{g_0}>> Y_0 = \P^n
\end{CD}
$$

Note that there is a natural action of $\PGL_{n+1}(\F_q)$ on $\P^n$,
and the above construction is equivariant with respect
to $\PGL_{n+1}(\F_q)$-action.
Hence we have a natural action of $\PGL_{n+1}(\F_q)$ on $B^n$.

\subsection{Divisors on $B^n$}
\label{SubsectionDivisors}

For $V \in \Grass_k(\P^n)(\F_q)$, we define a smooth irreducible
divisor $D_V$ on $B^n$ as follows.
If $k=n-1$, let $D_V \subset B^n$
be the strict transform of $V$.
If $k < n-1$, the strict transform $\widetilde{V} \subset Y_k$ of $V$
is a connected component of $Z_k$.
Therefore, $g_k^{-1}(\widetilde{V}) \subset Y_{k+1}$ is
a $\P^{n-d-1}$-bundle over $\widetilde{V}$.
Let $D_V \subset B^n$ be the strict transform of
$g_k^{-1}(\widetilde{V})$.

By construction, (\ref{CohomologyBlowUp3}) and induction,
an $\R$-divisor $D$ on $B^n$ is written {\em uniquely} as
$$ D = f^{\ast} \O(H) + \sum_{\substack{V \in \Grass_k(\P^n)(\F_q),\\
     0 \leq k \leq n-2}} a_V D_V
     \qquad \text{in} \quad H^2(B^n), $$
where $H$ is an $\R$-divisor on $\P^n$, and $a_V \in \R$.
We also see that,
all $\R$-divisors $D'$ on $B^n$ is written as
$$ D' = \sum_{V \in \Grass_{\ast}(\P^n)(\F_q)} a'_V D_V \qquad (a'_V \in \R)
   \qquad \text{in} \quad H^2(B^n). $$
Note that this second expression is {\em not unique}.

For $0 \leq k \leq n-1$,
we define a divisor $D_k$ on $B^n$ by
$$ D_k := \sum_{V \in \Grass_k(\P^n)(\F_q)} D_V. $$
By construction, we see that the support of $D_k$
is a disjoint union of smooth divisors on $B^n$.

In this section, we use the following terminology
for $\R$-divisors on $B^n$.

\begin{definition}
\label{DefinitionPositivePGLinvariant}
Let $D$ be an $\R$-divisor on $B^n$.
\begin{enumerate}
\item $D$ is called {\em positive}
  if $D$ can be written as
  $$ D = \sum_{V \in \Grass_{\ast}(\P^n)(\F_q)} a_V D_V
         \quad (a_V \in \R_{>0})
         \qquad \text{in} \quad H^2(B^n). $$
\item $D$ is called {\em $\PGL_{n+1}(\F_q)$-invariant}
  if the class of $D$ in $H^2(B^n)$ is
  $\PGL_{n+1}(\F_q)$-invariant.
\end{enumerate}
\end{definition}

It is easy to see that
$D$ is $\PGL_{n+1}(\F_q)$-invariant if and only if
$D$ can be written as
$ D = f^{\ast} \O(H) + \sum_{k=0}^{n-1} b_k D_k$
in $H^2(B^n)$ for some $H, b_k$
(for a similar criterion of positive divisors,
 see Proposition \ref{PGLinvariantProp}).

\subsection{Intersections of divisors $D_V$ on $B^n$}
\label{SubsectionCombinatorics}

Here we study the combinatorial structure of intersections
of divisors $D_V$ on $B^n$
(see also \cite{Mustafin}, proof of Theorem 4.1, I).

\begin{proposition}
\label{PropCombinatorics}
\begin{enumerate}
\item For $V,W \in \Grass_{\ast}(\P^n)(\F_q)$, $D_V \cap D_W \neq \emptyset$
  if and only if $V \subset W$ or $W \subset V$.
\item For $V \in \Grass_d(\P^n)(\F_q)$, we have a noncanonical isomorphism
  $$ D_V \cong B^d \times B^{n-d-1}, $$
  where $B^d$ (resp. $B^{n-d-1}$) is a variety
  constructed by successive blowing-ups
  of $\P^d$ (resp. $\P^{n-d-1}$) by the same way as $B^n$.
\item Let $x \in V$ be a point which does not lie in
  any linear subvariety defined over $\F_q$
  strictly contained in $V$.
  Let $\widetilde{V} \subset Y_d$ be the strict transform of $V$.
  Then $g_d \colon g_d^{-1}(\widetilde{V}) \to \widetilde{V}$
  is a trivial $\P^{n-d-1}$-bundle.
  Moreover, $\widetilde{V}$ is isomorphic to $B^d$,
  and an isomorphism $D_V \cong B^d \times B^{n-d-1}$ as in 2.
  is induced from a choice of $x \in V$ and a trivialization
    $g_d^{-1}(\widetilde{V})
        \cong \widetilde{V} \times \P(\check{N}_{{V/\P^n},x})$,
  where $\check{N}_{{V/\P^n},x}$ is the fiber at $x$
  of the conormal bundle $\check{N}_{V/\P^n}$.
\item For $W \in \Grass_{\ast}(\P^n)(\F_q)$ strictly contained in $V$,
  $D_W$ intersects transversally with $D_V$,
  and the intersection $D_W \cap D_V$ is of the form
  $$ D_W \cap D_V = \widetilde{D}_{W} \times B^{n-d-1}
     \qquad \text{on} \quad D_V, $$
  where $\widetilde{D}_{W}$ is a divisor on $B^d$
  corresponding to the inclusion $W \hookrightarrow V \cong \P^{d}$.
\item For $W' \in \Grass_{\ast}(\P^n)(\F_q)$ strictly containing $V$,
  $D_{W'}$ intersects transversally with $D_V$,
  and the intersection $D_{W'} \cap D_V$ is of the form
  $$ D_{W'} \cap D_V = B^d \times \widetilde{\widetilde{D}}_{W'}
     \qquad \text{on} \quad D_V, $$
  where $\widetilde{\widetilde{D}}_{W'}$ is
  a divisor on $B^{n-d-1}$ corresponding to the inclusion
  $\P(\check{N}_{{V/W'},x}) \hookrightarrow \P(\check{N}_{{V/\P^n},x})
   \cong \P^{n-d-1}$.
\end{enumerate}
\end{proposition}

\begin{proof}
Firstly, we shall prove the first assertion.
If $V \subset W$ or $W \subset V$, it is clear that
$D_V \cap D_W \neq \emptyset$.
Assume that $V \cap W$ is strictly smaller than $V, W$.
Then, after blowing-up along the strict transform of $V \cap W$
in the construction of $B^n$,
the strict transforms of $V$ and $W$ become disjoint.
Therefore, $D_V$ does not intersect with $D_W$.

Recall that, for a vector bundle $E$ of rank $k$
over a smooth variety $X$, $\P(E)$ is a trivial $\P^{k-1}$-bundle
if and only if $E \otimes_{\O_X} \L$ is
a trivial vector bundle for a line bundle $\L$ over $X$
(\cite{HartshorneAG}, II, Exercise 7.10).
We call such $E$ a {\em twist of a trivial vector bundle}
over $X$.
For $V \in \Grass_d(\P^n)(\F_q)$,
the normal bundle $N_{V/{\P^n}} \cong \O_V(1)^{\oplus n-d}$
is a twist of a trivial vector bundle over $V$.
Hence the normal bundle $N_{\widetilde{V}/Y_d}$
is also a twist of a trivial vector bundle over $Y_d$,
where $\widetilde{V} \subset Y_d$ is the strict transform of $V$
(\cite{Fulton}, B.6.10).
Therefore, $g_d^{-1}(\widetilde{V})$ is a trivial $\P^{n-d-1}$-bundle
over $\widetilde{V}$.

Here $\widetilde{V}$ is isomorphic to $B^d$.
This follows from the fact that, for a sequence of
regular embeddings $Z \hookrightarrow Y \hookrightarrow X$,
the strict transform of $Y$ in the blow-up of $X$ along $Z$
is isomorphic to the blow-up of $Y$ along $Z$
(\cite{Fulton}, B.6.9).

We fix a point $x \in V$ which does not lie in
any linear subvariety defined over $\F_q$
strictly contained in $V$,
and a trivialization
$ g_d^{-1}(\widetilde{V})
     \cong \widetilde{V} \times \P(\check{N}_{{V/\P^n},x}), $
where $\check{N}_{{V/\P^n},x}$ denotes the fiber at $x$
of the conormal bundle $\check{N}_{V/\P^n}$.
For $W \in \Grass_{\ast}(\P^n)(\F_q)$ strictly containing $V$,
a natural inclusion $V \subset W \subset \P^n$
induces $N_{V/W} \subset N_{V/\P^n}$ and hence
$\P(\check{N}_{{V/W}}) \subset \P(\check{N}_{{V/\P^n}})$.
Thus, $\P(\check{N}_{{V/W},x})$ is a linear subvariety
of $\P(\check{N}_{{V/\P^n},x})$.
From this, we see that the blowing-up of 
$g_d^{-1}(\widetilde{V})$ along the strict transform
of a linear subvariety strictly containing $V$
corresponds to the blowing-up
along the strict transform of
a linear subvariety of the second factor
$\P(\check{N}_{{V/\P^n},x})$ of
$g_d^{-1}(\widetilde{V}) \cong
\widetilde{V} \times \P(\check{N}_{{V/\P^n},x})$.
Therefore, we have an isomorphism
$D_V \cong B^d \times B^{n-d-1}$.
The remaining assertions follow from the construction
of this isomorphism.
\end{proof}

\begin{corollary}
\label{CorCombinatorics}
Let $H$ be an $\R$-divisor on $\P^n$,
$V \in \Grass_d(\P^n)(\F_q)$,
and $f \colon B^n \to \P^n$ a natural map
as in \S \ref{SubsectionConstruction}.
Fix an isomorphism $D_V \cong B^d \times B^{n-d-1}$
as in Proposition \ref{PropCombinatorics}.
Let $\pr_1 \colon D_V \to B^d$ be the projection
to the first factor.
Let $f' \colon B^d \to \P^d \cong V$ be
a map defined similarly as $f$.
Then we have
$$ \big( f^{\ast} \O(H) \big)\big|_{D_V}
     = \pr_1^{\ast} \Big( {f'}^{\ast} \big( \O(H)|_V \big) \Big)
   \qquad \text{in}\quad H^{2}(D_V). $$
\end{corollary}

\begin{proof}
This follows from the construction of
an isomorphism $D_V \cong B^d \times B^{n-d-1}$
in Proposition \ref{PropCombinatorics}.
\end{proof}

\subsection{$\PGL_{n+1}(\F_q)$-invariant $\R$-divisors on $B^n$}

\begin{proposition}
\label{PGLinvariantProp}
\begin{enumerate}
\item For $V \in \Grass_{n-1}(\P^n)(\F_q)$,
  we have
  $$ f^{\ast} \O(V) = \sum_{\substack{W \in \Grass_{\ast}(\P^n)(\F_q),
       \\ W \subset V}} D_W
       \qquad \text{in}\quad H^2(B^n). $$
\item For an $\R$-divisor $H$ on $\P^n$,
  $f^{\ast} \O(H)$ is $\PGL_{n+1}(\F_q)$-invariant.
\item Let $D$ be a $\PGL_{n+1}(\F_q)$-invariant
  $\R$-divisor on $B^n$ which is written as
  $$ D = \alpha\,f^{\ast} \O_{\P^n}(1) + \sum_{d=0}^{n-1} a_d D_d
     \quad (\alpha, a_d \in \R)
     \qquad \text{in}\quad H^2(B^n). $$
  Then, $D$ is positive
  if and only if $\alpha > 0$ and
  $$ a_d + \alpha\,\frac{|\P^{n-d-1}(\F_q)|}{|\P^n(\F_q)|} > 0
     \qquad \text{for all}\quad 0 \leq d \leq n-1. $$
\item Let $D$ be a $\PGL_{n+1}(\F_q)$-invariant $\R$-divisor
  on $B^n$. Take $V \in \Grass_d(\P^n)(\F_q)$,
  and fix an isomorphism $D_V \cong B^d \times B^{n-d-1}$
  as in Proposition \ref{PropCombinatorics}.
  Let $\pr_1 \colon D_V \to B^d,\ \pr_2 \colon D_V \to B^{n-d-1}$
  be projections.
  Then, there exist a $\PGL_{d+1}(\F_q)$-invariant
  $\R$-divisor $D'$ on $B^d$
  and a $\PGL_{n-d}(\F_q)$-invariant $\R$-divisor
  $D''$ on $B^{n-d-1}$ such that
  $$ \O(D)|_{D_V} = \pr_1^{\ast} \O(D') + \pr_2^{\ast} \O(D'')
     \qquad \text{in}\quad H^2(D_V). $$
\end{enumerate}
\end{proposition}

\begin{proof}
The first assertion follows from the fact that
the multiplicity of $V$ along $W \in \Grass_{\ast}(\P^n)(\F_q)$
with $W \subset V$ is equal to 1.

The second assertion is obvious because
$\PGL_{n+1}(\F_q)$ acts trivially on $H^2(\P^n)$
and $f^{\ast}$ is $\PGL_{n+1}(\F_q)$-equivariant.

We shall compute $f^{\ast} \O(H)$ explicitly.
Since $H^2(\P^n)$ is generated by the class of
$\O_{\P^n}(1)$, we write
$H = \alpha\,\O_{\P^n}(1)$ in $H^2(\P^n)$
for some $\alpha \in \R$.
Since $|\Grass_{n-1}(\P^n)(\F_q)| = |\P^n(\F_q)|$, and,
for each $W \in \Grass_d(\P^n)(\F_q)$,
there exist $|\P^{n-d-1}(\F_q)|$ elements in
$\Grass_{n-1}(\P^n)(\F_q)$ containing $W$,
we have the following equality in $H^2(B^n)$:
\begin{align*}
   f^{\ast} \O(H) &= \frac{\alpha}{|\P^n(\F_q)|}
     \sum_{V \in \Grass_{n-1}(\P^n)(\F_q)} \bigg(
     \sum_{\substack{W \in \Grass_{\ast}(\P^n)(\F_q),\\ W \subset V}}
       D_W \bigg) \\
   &= \frac{\alpha}{|\P^n(\F_q)|} \sum_{d=0}^{n-1} |\P^{n-d-1}(\F_q)|\,D_d.
\end{align*}
Of course, from this expression, we can check that
$f^{\ast} \O(H)$ is $\PGL_{n+1}(\F_q)$-invariant.

The third assertion easily follows from the above computation
because
$$ D = \sum_{d=0}^{n-1} \bigg( a_d +
     \alpha\,\frac{|\P^{n-d-1}(\F_q)|}{|\P^n(\F_q)|} \bigg) D_d
   \qquad \text{in}\quad H^2(B^n). $$

Finally, we shall prove the last assertion.
To avoid self-intersection, we write $D$ in the following form
$$ D = \alpha\,f^{\ast}\O_{\P^n}(1) +
         \sum_{0 \leq k \leq n-1,\ k \neq d} a_k D_k
   \quad (\alpha, a_k \in \R)
   \qquad \text{in}\quad H^2(B^n). $$
It is enough to treat the case $D = \alpha\,f^{\ast}\O_{\P^n}(1)$
and $D = D_k\ (k \neq d)$ separately.
The case $D = \alpha\,f^{\ast}\O_{\P^n}(1)$ follows from
Corollary \ref{CorCombinatorics}
and the second assertion.
The case $D = D_k\ (k \neq d)$ follows from
Proposition \ref{PropCombinatorics}.
\end{proof}

\subsection{Ample $\PGL_{n+1}(\F_q)$-invariant $\R$-divisors on $B^n$}

Here we study ample $\PGL_{n+1}(\F_q)$-invariant
$\R$-divisors on $B^n$.
First of all, we note that there exists at least
one ample $\PGL_{n+1}(\F_q)$-invariant $\R$-divisor on $B^n$.
It follows from the construction of $B^n$
in \S \ref{SubsectionConstruction} and
Proposition \ref{AmpleIntervalProp}.

\begin{proposition}
\label{AmplePGLinvariantProp}
Let $D$ be an ample $\PGL_{n+1}(\F_q)$-invariant $\R$-divisor on $B^n$.
Then, $D$ is positive
(see Definition \ref{DefinitionPositivePGLinvariant}).
\end{proposition}

\begin{proof}
We prove the assertion by induction on $n$.
The case $n=1$ is obvious.
Assume that Proposition \ref{AmplePGLinvariantProp} was
already proved in dimension $< n$.

Let $D$ be written uniquely as
$$ D = \alpha\,f^{\ast}\O_{\P^n}(1)
         + \sum_{d=0}^{n-2} a_d D_d
   \quad (\alpha, a_d \in \R)
   \qquad \text{in}\quad H^2(B^n). $$
By Proposition \ref{PGLinvariantProp},
$D$ is positive if and only if $\alpha > 0$ and
$$ a_d + \alpha\,\frac{|\P^{n-d-1}(\F_q)|}{|\P^n(\F_q)|} > 0
   \qquad \text{for all}\quad 0 \leq d \leq n-2. $$

Fix $V \in \Grass_{n-1}(\P^n)(\F_q)$.
By Proposition \ref{PropCombinatorics},
$D_V$ is isomorphic to $B^{n-1}$.
We consider the restriction of $\O(D)$ to $D_V$.
Let $f' \colon D_V \cong B^{n-1} \to \P^{n-1}$
be a map defined similarly as $f$.
Then, by 
Proposition \ref{PropCombinatorics} and
Corollary \ref{CorCombinatorics},
we have
$$ \O(D)|_{D_V} = \alpha\,{f'}^{\ast}\O_{\P^{n-1}}(1)
         + \sum_{d=0}^{n-2} a_d D'_d
   \qquad \text{in}\quad H^2(B^n), $$
where $D'_d$ is a divisor on $D_V \cong B^{n-1}$
defined similarly as $D_d$ on $B^n$.

By induction, $\O(D)|_{D_V}$ is positive.
Therefore, by Proposition \ref{PGLinvariantProp},
we have $\alpha > 0$ and
$$ a_d + \alpha\,\frac{|\P^{n-d-2}(\F_q)|}{|\P^{n-1}(\F_q)|} > 0
   \qquad \text{for all}\quad 0 \leq d \leq n-2. $$
Therefore, it is enough to show the following inequality
$$ \frac{|\P^{n-d-1}(\F_q)|}{|\P^{n}(\F_q)|}
     > \frac{|\P^{n-d-2}(\F_q)|}{|\P^{n-1}(\F_q)|}. $$
Since $|\P^k(\F_q)| = \frac{q^{k+1}-1}{q-1}$,
the above inequality is equivalent to
\begin{align*}
  &\hspace*{-0.2in}
        \frac{q^{n-d}-1}{q^{n+1}-1} > \frac{q^{n-d-1}-1}{q^n-1} \\
  &\iff q^{2n-d} - q^n - q^{n-d} + 1
           > q^{2n-d} - q^{n+1} - q^{n-d-1} + 1 \\
  &\iff q^{n+1} + q^{n-d-1} > q^n + q^{n-d} \\
  &\iff q^n(q-1) > q^{n-d-1}(q-1).
\end{align*}
Hence we prove Proposition \ref{AmplePGLinvariantProp}.
\end{proof}

\subsection{Hodge standard conjecture for $B^n$ with
ample $\PGL_{n+1}(\F_q)$-invariant $\R$-divisors}

\begin{lemma}
\label{LemmaRestriction}
For $0 \leq k \leq n-1$ and $a \in H^{k}(B^n)$,
if the restriction of $a$ to $D_V$ is zero
for all $V \in \Grass_{\ast}(\P^n)(\F_q)$, then $a=0$.
\end{lemma}

\begin{proof}
By applying Proposition \ref{BlowUpRestrictionProp}
successively, we see that $a$ is of the form
$a = f^{\ast} a'$ for some $a' \in H^{k}(\P^n)$.
For $V \in \Grass_{n-1}(\P^n)(\F_q)$,
the restriction of $a'$ to $V$ is zero because
the restriction of $a$ to $D_V$ is zero.
Since restriction map
$H^{k}(\P^n) \to H^k(V)$ is an isomorphism,
we conclude that $a=0$.
\end{proof}

The following proposition is the main result of this section.

\begin{proposition}
\label{MainPropLefschetzHodge}
Let $D$ be an ample $\PGL_{n+1}(\F_q)$-invariant $\R$-divisor on $B^n$.
Then the hard Lefschetz conjecture
(Conjecture \ref{HardLefschetzConjecture})
and the Hodge standard conjecture
(Conjecture \ref{HodgeConj}) hold for $(B^n,D)$.
\end{proposition}

\begin{proof}
We prove the assertion by induction on $n$.
The case $n=1$ is obvious.
Assume that Proposition \ref{MainPropLefschetzHodge} was
already proved in dimension $<n$.

By the construction of $B^n$
in \S \ref{SubsectionConstruction} and
Proposition \ref{BlowUpHodgeConjProp},
we see that there exists at least
one ample $\PGL_{n+1}(\F_q)$-invariant $\R$-divisor on $B^n$
for which Proposition \ref{MainPropLefschetzHodge} holds.
If $D, D'$ are ample $\PGL_{n+1}(\F_q)$-invariant $\R$-divisors on $B^n$,
then $t D + (1-t) D'$ is also
an ample $\PGL_{n+1}(\F_q)$-invariant $\R$-divisor on $B^n$
for $t \in \R$ with $0 \leq t \leq 1$.
Therefore, by Corollary \ref{IndependenceCor},
it is enough to show that
the hard Lefschetz conjecture
(Conjecture \ref{HardLefschetzConjecture})
holds for $(B^n,D)$ for
all ample $\PGL_{n+1}(\F_q)$-invariant $\R$-divisors $D$ on $B^n$.

Assume that there exists an ample $\PGL_{n+1}(\F_q)$-invariant
$\R$-divisor $L$ on $B^n$ and a nonzero cohomology class
$a \in H^{k}(B^n)$ for $0 \leq k \leq n-1$ such that
$$ L^{n-k} \cup a = 0. $$

For $0 \leq d \leq n-1$,
take $V \in \Grass_d(\P^n)(\F_q)$ and
fix an isomorphism $D_V \cong B^d \times B^{n-d-1}$
as in Proposition \ref{PropCombinatorics}.
Let $\pr_1 \colon D_V \to B^d,\ \pr_2 \colon D_V \to B^{n-d-1}$
be projections.
Then, by Proposition \ref{PGLinvariantProp},
the restriction $\O(L)|_{D_V}$ can be written as
$$ \O(L)|_{D_V} = \pr_1^{\ast} \O(D') + \pr_2^{\ast} \O(D'')
     \qquad \text{in}\quad H^2(D_V), $$
where $D'$ (resp. $D''$) is an ample
$\PGL_{d+1}(\F_q)$-invariant
(resp. $\PGL_{n-d}(\F_q)$-invariant)
$\R$-divisor on $B^d$ (resp. $B^{n-d-1}$).
Therefore, by induction hypothesis and
Proposition \ref{ConjectureProduct},
the hard Lefschetz conjecture
(Conjecture \ref{HardLefschetzConjecture})
and the Hodge standard conjecture
(Conjecture \ref{HodgeConj})
hold for $\big( D_V, \O(L)|_{D_V} \big)$.

By restricting $L^{n-k} \cup a = 0$ to $D_V$,
we have $\big( \O(L)|_{D_V} \big)^{n-k} \cup (a|_{D_V}) = 0$.
Since $\dim D_V = n-1$,
this implies $a|_{D_V} \in H^{k}(D_V)$ is
in the primitive part $P^{k}(D_V)$.
Therefore, by the Hodge standard conjecture
(Conjecture \ref{HodgeConj})
for $\big( D_V, \O(L)|_{D_V} \big)$,
we have
$$ (-1)^{k/2} \cdot \big( \O(L)|_{D_V} \big)^{n-k-1}
       \cup (a|_{D_V}) \cup (a|_{D_V}) \geq 0 $$
and the equality holds if and only if $a|_{D_V} = 0$.

By Proposition \ref{AmplePGLinvariantProp},
$L$ can be written as
$$ L = \sum_{V \in \Grass_{\ast}(\P^n)(\F_q)} a_V D_V
   \quad (a_V \in \R_{>0})
   \qquad \text{in} \quad H^2(B^n). $$
Therefore, we have
\begin{align*}
    &\hspace*{-0.3in} (-1)^{k/2} \cdot L^{n-k} \cup a \cup a \\
    &= \sum_{V \in \Grass_{\ast}(\P^n)(\F_q)} a_V \cdot D_V \cup
        \big( (-1)^{k/2} \cdot L^{n-k-1} \cup a \cup a \big) \\
    &= \sum_{V \in \Grass_{\ast}(\P^n)(\F_q)} a_V \Big(
        (-1)^{k/2} \cdot \big( \O(L)|_{D_V} \big)^{n-k-1}
          \cup (a|_{D_V}) \cup (a|_{D_V}) \Big) \\
    &= 0
\end{align*}
because $L^{n-k} \cup a = 0$.
Since $a_V>0$ for all $V \in \Grass_{\ast}(\P^n)(\F_q)$,
we have $a|_{D_V} = 0$ for all $V \in \Grass_{\ast}(\P^n)(\F_q)$.
By Lemma \ref{LemmaRestriction}, we have $a=0$.
This is a contradiction.
Hence the proof of Proposition \ref{MainPropLefschetzHodge}
is complete.
\end{proof}

\begin{remark}
\label{RemarkAmpleDivisorP-adicUniformization}
Mustafin, Kurihara showed that
$B^d$ is isomorphic to 
all irreducible components of the special fiber
of the formal scheme model $\widehat{\Omega}^d_{\O_K}$
over $\text{Spf}\,\O_K$
of the Drinfeld upper half space
$\widehat{\Omega}^d_{K}$ of dimension $d$
(for details, see \S \ref{SectionP-adic}, see also
 \cite{Mustafin}, \cite{Kurihara}).
Moreover, they also showed that
  $$ D = -(d+1) f^{\ast} \O_{\P^d}(1) + \sum_{k=0}^{d-1} (d-k) D_k $$
is an ample divisor on $B^d$, and
$D$ coincides with the restriction
of the \lq\lq relative dualizing sheaf''
$\omega_{\widehat{\Omega}^d_{\O_K}/\O_K}$
to an irreducible component of the special fiber.
These facts are crucial in our application
to $p$-adic uniformization in \S \ref{SectionP-adic}.
\end{remark}

\begin{remark}
In our proof of Proposition \ref{AmplePGLinvariantProp},
we heavily use the assumption that $D$ is
$\PGL_{n+1}(\F_q)$-invariant.
Therefore, our proof does not work for ample
non-$\PGL_{n+1}(\F_q)$-invariant $\R$-divisors on $B^n$.
It is easy to see that if all ample $\R$-divisors on
$B^n$ are positive,
then Proposition \ref{MainPropLefschetzHodge} holds for them.
However, the author does not even know whether
it is natural to expect this. He has neither evidence
nor counter-example for ample
non-$\PGL_{n+1}(\F_q)$-invariant $\R$-divisors.
\end{remark}

\section{Review of the weight spectral sequence of Rapoport-Zink}
\label{SectionFiltration}

In this section, we recall the definitions and
basic properties of monodromy filtration, weight filtration,
and the weight spectral sequence of Rapoport-Zink.
Since we work over a local field in this section,
we slightly change the notation.
As in \S \ref{SectionIntroduction},
let $K$ be a complete discrete valuation field
with finite residue field $\F_q$ of characteristic $p>0$,
$l$ be a prime number different from $p$.
Let $X$ be a proper smooth variety of dimension $n$
over $K$, and $V := H^w_{\text{\rm \'et}}(X_{\overline{K}},\Q_l)$.

\subsection{Monodromy filtration}
\label{SectionMonodromyFiltration}

Let $I_K$ be the inertia group of $K$,
which is a subgroup of $\Gal(\overline{K}/K)$
defined by the exact sequence:
$$
\begin{CD}
1 @>>> I_K @>>> \Gal(\overline{K}/K) @>>> \Gal(\overline{\F}_q/\F_q) @>>> 1.
\end{CD}
$$
$\Gal(\overline{\F}_q/\F_q)$ acts on $I_K$ by conjugation
($\tau \colon \sigma \mapsto \tau \sigma \tau^{-1},
\ \tau \in \Gal(\overline{\F}_q/\F_q),\ \sigma \in I_K$).
The pro-$l$-part of $I_K$ is isomorphic to $\Z_l(1)$
as a $\Gal(\overline{\F}_q/\F_q)$-module by
$$ t_l \colon I_K \ni \sigma \mapsto
   \left( \frac{\sigma(\pi^{1/l^m})}{\pi^{1/l^m}} \right)_{\!\!m}
   \in \varprojlim \mu_{l^m} = \Z_l(1), $$
where $\pi$ is a uniformizer of $K$,
and $\mu_{l^m}$ is the group of $l^m$-th roots of unity.
It is known that $t_l$ is independent of
the choice of $\pi$ and its $l^m$-th root $\pi^{1/l^m}$ (\cite{Serre1}).

By Grothendieck's monodromy theorem (\cite{SerreTate}, Appendix),
there exist $r,s \geq 1$ such that
$\big( \rho(\sigma)^r - 1 \big)^s = 0$
for all $\sigma \in I_K$.
Therefore, by replacing $K$ by its finite extension,
$I_K$ acts on $V$ through
$t_l \colon I_K \to \Z_l(1)$ and this action is unipotent.
Then there is a unique nilpotent map
of $\Gal(\overline{K}/K)$-representations
called the {\em monodromy operator} $N \colon V(1) \to V$
such that
$$ \rho(\sigma)
    = \exp \big( t_l(\sigma) N \big)
   \quad \text{for all} \quad \sigma \in I_K. $$
Here $N$ is a nilpotent map means that
$N^r \colon V(r) \to V$ is zero for some $r \geq 1$.

\begin{definition} (\cite{WeilII}, I, 1.7.2)
\label{DefMonodromyFiltration}
There exists a unique filtration $M_{\bullet}$ on $V$
called the {\em monodromy filtration} characterized by
the following properties.
\begin{enumerate}
\item $M_{\bullet}$ is an increasing filtration
  $\cdots \subset M_{i-1} V \subset M_i V \subset M_{i+1} V \subset \cdots$
  of $\Gal(\overline{K}/K)$-representations such that
  $M_i V = 0$ for sufficiently small $i$ and
  $M_i V = V$ for sufficiently large $i$.
\item $N \big( M_i V(1) \big) \subset M_{i-2} V$ for all $i$.
\item By the second condition,
  we can define $N \colon \Gr^M_{i} V (1) \to \Gr^M_{i-2} V$,
  where $\Gr^M_{i} V := M_{i} V / M_{i-1} V$.
  Then, for each $r \geq 0$,
  $N^r \colon \Gr^M_{r} V (r) \to \Gr^M_{-r} V$
  is an isomorphism.
\end{enumerate}
\end{definition}

\begin{remark}
\label{MonodromyFiltrationFiniteExtension}
In the above definition, we replace $K$ by its finite extension.
We can easily see that
$M_{\bullet}$ is stable under the action of $\Gal(\overline{K}/K)$
for the original $K$.
Therefore, we can define the monodromy filtration $M_{\bullet}$
as a filtration of $\Gal(\overline{K}/K)$-representations
without replacing $K$ by its finite extension.
\end{remark}

\subsection{Weight filtration}
\label{SectionWeightFiltration}

Let $\Fr_q \in \Gal(\overline{\F}_q/\F_q)$ be
the inverse of the $q$-th power map on $\overline{\F}_q$
called the {\em geometric Frobenius element}.
A $\Gal(\overline{\F}_q/\F_q)$-representation
is said to have {\em weight $k$} if all eigenvalues of
the action of $\Fr_q \in \Gal(\overline{\F}_q/\F_q)$
are algebraic integers whose all complex conjugates have
complex absolute value $q^{k/2}$.

\begin{definition} (\cite{HodgeI}, \cite{WeilII}, I, 1.7.5)
\label{DefWeightFiltration}
There exists a unique filtration $W_{\bullet}$ 
called the {\em weight filtration} on $V$ characterized by
the following properties
(for existence, see \S \ref{SubsectionWeightSpectralSequence}).
\begin{enumerate}
\item $W_{\bullet}$ is an increasing filtration
  $\cdots \subset W_{i-1} V \subset W_i V \subset W_{i+1} V \subset \cdots$
  of $\Gal(\overline{K}/K)$-representations such that
  $W_i V = 0$ for sufficiently small $i$ and
  $W_i V = V$ for sufficiently large $i$.
\item For a lift $\widetilde{\Fr}_q$ of $\Fr_q$ in $\Gal(\overline{K}/K)$,
  all eigenvalues of the action of $\widetilde{\Fr}_q$
  on each $\Gr^W_{i} V := W_i V / W_{i-1} V$ are algebraic integers
  whose all complex conjugates have complex absolute value $q^{i/2}$.
\end{enumerate}
\end{definition}

\begin{remark}
By the second condition, the monodromy operator
$N \colon \Gr^W_{i} V(1) \to \Gr^W_{i} V$ is zero for each $i$
(see \S \ref{SectionMonodromyFiltration}).
Hence, by replacing $K$ by its finite extension $K'$ with
residue field $\F_{q^r}$, $I_{K'}$ acts on $\Gr^W_{i} V$ trivially,
and the $\Gal(\overline{\F}_{q^r}/\F_{q^r})$-representation
$\Gr^W_{i} V$ has weight $i$ for each $i$.
\end{remark}

\subsection{Weight spectral sequence of Rapoport-Zink}
\label{SubsectionWeightSpectralSequence}

\begin{definition}
\label{DefSemistableModel}
A regular scheme $\X$ which is proper and flat over $\O_K$
is called a {\em proper semistable model}
of $X$ over $\O_K$ if the generic fiber
$\X_K := \X \otimes_{\O_K} K$ is isomorphic to $X$
and the special fiber
$\X_{\F_q} := \X \otimes_{\O_K} \F_q$
is a divisor of $\X$ with normal crossings.
Moreover, if $\X_{\F_q}$ is a divisor of $\X$
with simple normal crossings,
$\X$ is called a {\em proper strictly semistable model}
of $X$ over $\O_K$.
\end{definition}

We recall the weight spectral sequence of
Rapoport-Zink (\cite{Rapoport-Zink1}).
Assume that $X$ has a proper strictly semistable model $\X$ over $\O_K$.
Let $X_1,\ldots,X_m$ be the irreducible components of
the special fiber of $\X$, and
$$ X^{(k)} := \coprod_{1 \leq i_1 < \cdots < i_k \leq m}
     X_{i_1} \cap \cdots \cap X_{i_k}. $$
Then $X^{(k)}$ is a disjoint union of proper smooth irreducible varieties
of dimension $n-k+1$ over $F$.
The {\em weight spectral sequence of Rapoport-Zink}
is as follows:
\begin{center}
\begin{tabular}{l}
$\displaystyle E_1^{-r,\,w+r} = \bigoplus_{k \geq \max\{0,-r\}}
      H^{w-r-2k}_{\text{\rm \'et}}
      \big( X^{(2k+r+1)}_{\overline{\F}_q},\, \Q_l(-r-k) \big)$ \\
\hspace*{3.2in} $\displaystyle \Longrightarrow \quad
      H^w_{\text{\rm \'et}}(X_{\overline{K}},\Q_l)$.
\end{tabular}
\end{center}
This spectral sequence is $\Gal(\overline{K}/K)$-equivariant.
The map $d_1^{i,\,j} \colon E_1^{i,\,j} \to E_1^{i+1,\,j}$
can be described in terms of restriction morphisms and Gysin morphisms
explicitly (see \cite{Rapoport-Zink1}, 2.10 for details).

The action of the monodromy operator $N$ on
$H^w_{\text{\rm \'et}}(X_{\overline{K}},\Q_l)$
in \S \ref{SectionMonodromyFiltration}
is induced from a natural map
$N \colon E_1^{i,\,j}(1) \to E_1^{i+2,\,j-2}$
satisfying
$$
\begin{CD}
  N^r \colon E_1^{-r,\,w+r}(r) @>{\cong}>> E_1^{r,\,w-r}
\end{CD}
$$
for all $r,w$.
We can describe $N \colon E_1^{i,\,j}(1) \to E_1^{i+2,\,j-2}$
explicitly (\cite{Rapoport-Zink1}, 2.10).
Hence we can define $N$ without replacing $K$ by
its finite extension
(see Remark \ref{MonodromyFiltrationFiniteExtension}).

The inertia group
$I_K$ acts on each $E_1^{i,\,j}$ trivially and
$\Gal(\overline{\F}_q/\F_q)$ acts on them.
By the Weil conjecture (\cite{WeilI}, \cite{WeilII}),
$H^{w-r-2k}_{\text{\rm \'et}}
  \big( X^{(2k+r+1)}_{\overline{\F}_q},\, \Q_l(-r-k) \big)$
has weight $(w-r-2k) -2(-r-k) = w+r$
(see \S \ref{SectionWeightFiltration}).
Hence $E_1^{i,\,j}$ has weight $j$.
Therefore, the filtration on $H^w_{\text{\rm \'et}}(X_{\overline{K}},\Q_l)$
induced by the weight spectral sequence is the weight filtration
$W_{\bullet}$ in Definition \ref{DefWeightFiltration}.
This proves the existence of $W_{\bullet}$
in Definition \ref{DefWeightFiltration}
(for general $X$, we may use de Jong's alteration \cite{deJong}).
Moreover, since
$d_r^{i,\,j} \colon E_r^{i,\,j} \to E_r^{i+r,\,j-r+1}$
is a map between $\Gal(\overline{\F}_q/\F_q)$-representations
with different weights for $r \geq 2$,
$d_r^{i,\,j} = 0$ for $r \geq 2$.
Therefore, the weight spectral sequence degenerates at $E_2$.

By combining above facts, we see that
the weight-monodromy conjecture (Conjecture \ref{WMC})
is equivalent to the following conjecture
on the weight spectral sequence of Rapoport-Zink.

\begin{conjecture}[\cite{Rapoport-Zink1}, \cite{Illusie1}, \cite{Illusie2}, \cite{Illusie3}]
\label{WMC_semistable}
Let $X$ be a proper smooth variety over $K$
which has a proper strictly semistable model $\X$ over $\O_K$.
Let
$E_1^{-r,\,w+r} \Rightarrow
    H^w_{\text{\rm \'et}}(X_{\overline{K}},\Q_l)$
be the weight spectral sequence of Rapoport-Zink.
Then $N^r$ induces an isomorphism
$$
\begin{CD}
  N^r \colon E_2^{-r,\,w+r}(r) @>{\cong}>> E_2^{r,\,w-r}
\end{CD}
$$
on $E_2$-terms for all $r,w$.
\end{conjecture}

\subsection{A positivity argument of Steenbrink, M. Saito}
\label{SubsectionSteenbrinkSaito}

\begin{proposition}
\label{WMCHodgeProp}
Let $X$ be a proper smooth variety over $K$
which has a proper strictly semistable model $\X$ over $\O_K$.
Let $X_1,\ldots,X_m$ be the irreducible components of
the special fiber of $\X$.
Assume that $\X$ is projective over $\O_K$ with
an ample line bundle $\L$.
Assume further that, for $1 \leq i_1 < \cdots < i_k \leq m$,
each irreducible component $Y$ of $X_{i_1} \cap \cdots \cap X_{i_k}$
satisfies Assumption \ref{AssumptionCohomology}
in \S \ref{SectionNotation},
and the Hodge standard conjecture
(Conjecture \ref{HodgeConj}) holds for $(Y,\L|_{Y})$.
Then, Conjecture \ref{WMC_semistable} holds for $\X$,
and Conjecture \ref{WMC} holds for $X$.
\end{proposition}

In the followings,
we reduce the weight-monodromy conjecture
(Conjecture \ref{WMC}, Conjecture \ref{WMC_semistable})
to an assertion that the restriction of a nondegenerate pairing
to a subspace is nondegenerate
(see Lemma \ref{LemmaNondegeneracy} below).
This method has a long story starting in SGA7,
where Grothendieck considered the case of curves
and abelian varieties (\cite{SGA7-I}, see also \cite{Illusie1}).

The outline of the proof here is essentially the same as in
\cite{MSaito}, 4.2.5, where M. Saito used polarized
Hodge structures to prove a Hodge analogue of
Conjecture \ref{WMC_semistable}
(see also \cite{MSaito2}).
The only difference is that we use $\Q$-structures
on $l$-adic cohomology and the Hodge standard
conjecture instead of the polarized Hodge structures
(A similar argument can be found in \cite{Ito2}).
Since the presentation of M. Saito in \cite{MSaito}
seems too sophisticated, it is not very clear to non-specialists
that the argument in \cite{MSaito} can also
be applied to the situation in Proposition \ref{WMCHodgeProp}
even in positive or mixed characteristic.
So we reproduce the argument here with slight modification
suitable for Proposition \ref{WMCHodgeProp}
for the reader's convenience.
Therefore, those who are familiar with the argument
in \cite{MSaito} may skip to the next section.

First of all, we recall the structure of the $E_1$-terms of
the weight spectral sequence of Rapoport-Zink.
Let
$$
\begin{CD}
\displaystyle
d_1^{-r,\,w+r} = \sum_{k \geq \max\{0,-r\}} (\rho + \tau) \colon
  E_1^{-r,\,w+r} @>>> E_1^{-r+1,\,w+r}
\end{CD}
$$
denote the differential on $E_1$-terms, where
$$
\begin{CD}
\rho \colon H^s_{\text{\rm \'et}}(X^{(t)}_{\overline{\F}_q},\Q_l) @>>>
  H^s_{\text{\rm \'et}}(X^{(t+1)}_{\overline{\F}_q},\Q_l)
\end{CD}
$$
is a linear combination of restriction morphisms for some $s,t$
($\rho$ is $(-1)^{r+k} \theta$ in \cite{Rapoport-Zink1}, 2.10),
and
$$
\begin{CD}
\tau \colon H^s_{\text{\rm \'et}}(X^{(t)}_{\overline{\F}_q},\Q_l) @>>>
  H^{s+2}_{\text{\rm \'et}}(X^{(t-1)}_{\overline{\F}_q},\Q_l(1))
\end{CD}
$$
is a linear combination of Gysin morphisms for some $s,t$
($\tau$ is $(-1)^k d_1'$ in \cite{Rapoport-Zink1}, 2.10).
From the construction,
$\rho$ increases the index $k$ by 1,
and $\tau$ preserves the index $k$.
Moreover, $\rho, \tau$ satisfy
$\rho \circ \rho = 0,\ \tau \circ \tau = 0,
\ \tau \circ \rho + \rho \circ \tau = 0$.
Hence $\rho, \tau$ also satisfy
$$ \rho \circ \tau \circ \rho = - \rho \circ \rho \circ \tau = 0,
   \qquad
   \tau \circ \rho \circ \tau = - \tau \circ \tau \circ \rho = 0. $$

Let $k$ be an integer with $1 \leq k \leq n$.
To distinguish $\rho$, $\tau$ for different degree,
we use the following notation.
$$
\begin{CD}
   \rho_i^{(k)} \colon
        H^i_{\text{\rm \'et}}(X^{(k)}_{\overline{\F}_q},\Q_l) @>>>
        H^i_{\text{\rm \'et}}(X^{(k+1)}_{\overline{\F}_q},\Q_l), \\
   \tau_i^{(k+1)} \colon
        H^i_{\text{\rm \'et}}(X^{(k+1)}_{\overline{\F}_q},\Q_l) @>>>
        H^{i+2}_{\text{\rm \'et}}(X^{(k)}_{\overline{\F}_q},\Q_l(1)).
\end{CD}
$$
In the followings, to simply the notation,
we fix an isomorphism $\Q_l \cong \Q_l(1)$
to ignore the Tate twists, and simply write
$H^{i}_{\text{\rm \'et}}(X^{(k)})
  := H^{i}_{\text{\rm \'et}}(X^{(k)}_{\overline{\F}_q}, \Q_l)$,
where we add a subscript \lq\lq \'et''
not to confuse with the notation $H^i(X)$
in \S \ref{SectionNotation}.

Let
$L \colon H^i_{\text{\rm \'et}}(X^{(k)})
\to H^{i+2}_{\text{\rm \'et}}(X^{(k)})$ be the Lefschetz operator
defined by $\L$.
Since $L$ commutes with $\rho, \tau$, we use the same letter
$L$ for different $i, k$, which does not cause any confusion.
Let
$$
  H^i_{\text{\rm \'et}}(X^{(k)})
    \times H^{2(n-k+1)-i}_{\text{\rm \'et}}(X^{(k)})
     \to \Q_l, \quad
  H^i_{\text{\rm \'et}}(X^{(k+1)})
    \times H^{2(n-k)-i}_{\text{\rm \'et}}(X^{(k+1)})
     \to \Q_l
$$
be the sum of the cup product pairings.
For simplicity, we denote them by $\cup$.
Then, for
$a \in H^i_{\text{\rm \'et}}(X^{(k)})$,
\ $b \in H^{2(n-k)-i}_{\text{\rm \'et}}(X^{(k+1)})$,
we have
\begin{equation}
\label{CupProductDuality}
  a \cup \tau_{2(n-k)-i}^{(k+1)}(b)
       = \pm \rho_i^{(k)}(a) \cup b.
\end{equation}
Namely, $\rho_i^{(k)}$ and $\tau_{2(n-k)-i}^{(k+1)}$
are dual to each other with respect to $\cup$ up to sign.
Here we do not specify the sign $\pm$ because it does not matter
in the following proof.

\begin{lemma}
\label{WeightSpectralSequenceLemma}
In the following sequence of $\rho, \tau$
$$
\begin{CD}
  H^{i-2}_{\text{\rm \'et}}(X^{(k+1)}) @>{\tau_{i-2}^{(k+1)}}>>
  H^i_{\text{\rm \'et}}(X^{(k)}) @>{\rho_i^{(k)}}>>
  H^i_{\text{\rm \'et}}(X^{(k+1)}) @>{\tau_i^{(k+1)}}>>
  H^{i+2}_{\text{\rm \'et}}(X^{(k)}),
\end{CD}
$$
we have
  $ \big( \Ker \tau_i^{(k+1)} \big) \cap \big( \Im \rho_i^{(k)} \big)
        = \Im \big( \rho_i^{(k)} \circ \tau_{i-2}^{(k+1)} \big). $
Moreover, in the following sequence of $\rho, \tau$
$$
\begin{CD}
  H^i_{\text{\rm \'et}}(X^{(k)}) @>{\rho_i^{(k)}}>>
  H^i_{\text{\rm \'et}}(X^{(k+1)}) @>{\tau_i^{(k+1)}}>>
  H^{i+2}_{\text{\rm \'et}}(X^{(k)}) @>{\rho_{i+2}^{(k)}}>>
  H^{i+2}_{\text{\rm \'et}}(X^{(k+1)}),
\end{CD}
$$
we have
$\big( \Ker \rho_{i+2}^{(k)} \big) \cap \big( \Im \tau_i^{(k+1)} \big)
    = \Im \big( \tau_i^{(k+1)} \circ \rho_i^{(k)} \big)$.
\end{lemma}

By assuming the above lemma, it is easy to prove
Proposition \ref{WMCHodgeProp}.

\begin{proof}[Proof of Proposition \ref{WMCHodgeProp}
assuming Lemma \ref{WeightSpectralSequenceLemma}]
For $r \geq 1$, we write
$E_1^{-r,\,w+r},\ E_1^{r,\,w-r}$ explicitly
as follows
(see \S \ref{SubsectionWeightSpectralSequence})
$$ E_1^{-r,\,w+r} = \bigoplus_{k \geq 0}
   H^{w-r-2k}_{\text{\rm \'et}}(X^{(2k+r+1)}), \qquad
   E_1^{r,\,w-r} = \bigoplus_{k \geq r}
   H^{w+r-2k}_{\text{\rm \'et}}(X^{(2k-r+1)}). $$
Put $k' = k-r$ and rewrite $E_1^{r,\,w-r}$ as
$$ E_1^{r,\,w-r} = \bigoplus_{k' \geq 0}
     H^{w-r-2k'}_{\text{\rm \'et}}(X^{(2k'+r+1)}). $$
From this expression, there is a natural
isomorphism between $E_1^{-r,\,w+r}$ and $E_1^{r,\,w-r}$,
which is nothing but an isomorphism induced by $N^r$
in \S \ref{SubsectionWeightSpectralSequence}.
Now, we look at the following part of
the weight spectral sequence of Rapoport-Zink:
\begin{equation}
\label{WeightSpectralSequencePart}
\xymatrix{
    E_1^{-r-1,\,w+r} \ar[r] \ar[dr]^{N^r} &
    E_1^{-r,\,w+r} \ar[r] \ar[dr]^{N^r} &
    E_1^{-r+1,\,w+r} \ar[dr]^{N^r} \\
  & E_1^{r-1,\,w-r} \ar[r] &
    E_1^{r,\,w-r}   \ar[r] &
    E_1^{r+1,\,w-r}.
}
\end{equation}
Note that $N^r$ is a morphism of complexes from the upper line
to the lower one.
It is easy to see that the leftmost $N^r$ is surjective,
the middle $N^r$ is an isomorphism,
and the rightmost $N^r$ is injective.
To simplify the notation, we define
\begin{align*}
  A &:= \bigoplus_{k \geq 0} H^{w-r-2k-2}_{\text{\rm \'et}}(X^{(2k+r+2)}), \\
  B &:= \bigoplus_{k \geq 1} H^{w-r-2k}_{\text{\rm \'et}}(X^{(2k+r+1)}),
  \quad
  B' := H^{w-r}_{\text{\rm \'et}}(X^{(r+1)}), \\
  C &:= \bigoplus_{k \geq 1} H^{w-r-2k+2}_{\text{\rm \'et}}(X^{(2k+r)}),
  \quad
  D := H^{w-r+2}_{\text{\rm \'et}}(X^{(r)}),
  \quad
  E := H^{w-r}_{\text{\rm \'et}}(X^{(r)}).
\end{align*}
Then, (\ref{WeightSpectralSequencePart}) becomes
\begin{equation}
\label{WeightSpectralSequencePart2}
\xymatrix{
    A \ar[rr]^-{d+d'} \ar[drr] & &
    B \oplus B' \ar[rr]^-{(e+0,\,f+f')} \ar[drr] & &
    C \oplus D \ar[drr] \\
& & A \oplus E \ar[rr]_-{(d+d',\,0+g)} & &
    B \oplus B' \ar[rr]_-{(e, f)} & &
    C
}
\end{equation}
for some
$d \colon A \to B,\ d' \colon A \to B',\ e \colon B \to C,
\ f \colon B' \to C,\ f' \colon B' \to D,\ g \colon E \to B'$,
where $B \oplus B' \to B \oplus B'$ is the identity map,
and $A \to A \oplus E$ (resp. $C \oplus D \to C$) is
a natural map $a \mapsto (a,0)$
(resp. $(c,d) \mapsto c$).

It is enough to show that the identity map
$B \oplus B' \to B \oplus B'$ induces an isomorphism
between the cohomologies of the first and second rows
of (\ref{WeightSpectralSequencePart2}).

We shall prove the injectivity.
Take an element $(b,b') \in B \oplus B'$
such that $e(b) + f(b') = 0,\ f'(b') = 0$.
Assume that $(b + b') \in \Im(d+d', 0+g)$.
Namely, there exists $(x,y) \in A \oplus E$
such that $d(x) = b,\ d'(x) + g(y) = b'$.
Then we have
\begin{align*}
   e(b) + f(b') &= e \circ d (x) + f \circ d'(x)
                     + f \circ g(y) = 0 \\
   f'(b') &= f' \circ d'(x) + f' \circ g(y) = 0.
\end{align*}
From the definition, we see that
$d' = \tau_{w-r-2}^{(r+2)},\ f' = \tau_{w-r}^{(r+1)},
\ g = \rho_{w-r}^{(r)}$.
Hence
$f' \circ d' = \tau_{w-r}^{(r+1)} \circ \tau_{w-r-2}^{(r+2)} = 0$.
Therefore, we have
$$ f' \circ g(y)
    = \tau_{w-r}^{(r+1)} \circ \rho_{w-r}^{(r)}(y) = 0. $$
By the first assertion of Lemma \ref{WeightSpectralSequenceLemma},
there exists $z \in H^{w-r-2}_{\text{\rm \'et}}(X^{(r+1)})$
such that
$$ g(y) = \rho_{w-r}^{(r)}(y)
      = \rho_{w-r}^{(r)} \circ \tau_{w-r-2}^{(r+1)}(z). $$
Since
$\rho_{w-r}^{(r)} \circ \tau_{w-r-2}^{(r+1)}(z)
  = - \tau_{w-r-2}^{(r+2)} \circ \rho_{w-r-2}^{(r+1)}(z)$,
we have
$$ b' = d'(x) + g(y) = \tau_{w-r-2}^{(r+2)}(x)
    - \tau_{w-r-2}^{(r+2)} \circ \rho_{w-r-2}^{(r+1)}(z)
      = \tau_{w-r-2}^{(r+2)}\big(x - \rho_{w-r-2}^{(r+1)}(z)\big). $$
On the other hand,
the restriction of $d \colon A \to B$ to
the $k=0$ part $H^{w-r-2}_{\text{\rm \'et}}(X^{(r+2)})$
is $\rho_{w-r-2}^{(r+2)}$.
Then, we have
$$ d\,\big(x - \rho_{w-r-2}^{(r+1)}(z)\big)
  = d(x) - \rho_{w-r-2}^{(r+2)} \circ \rho_{w-r-2}^{(r+1)}(z)
  = d(x) = b, $$
since $\rho_{w-r-2}^{(r+2)} \circ \rho_{w-r-2}^{(r+1)} = 0$.
Hence $(b, b')$ is the image of
$x - \rho_{w-r-2}^{(r+1)}(z)$ by $d+d'$,
i.e. $(b, b') \in \Im(d+d')$.
This proves the injectivity.

Finally, we shall prove the surjectivity.
Take an element $(b,b') \in B \oplus B'$
with $e(b) + f(b') = 0$.
From the construction, we see that $f = \rho_{w-r}^{(r+1)}$
and the restriction of $e \colon B \to C$ to the $k=1$ part
$H^{w-r-2}_{\text{\rm \'et}}(X^{(r+3)})$ is $\tau_{w-r-2}^{(r+3)}$.
Therefore, the equality $e(b) + f(b') = 0$
implies
$\tau_{w-r-2}^{(r+3)}(\tilde{b}) + \rho_{w-r}^{(r+1)}(b') = 0$
in the $k=1$ part $H^{w-r}_{\text{\rm \'et}}(X^{(r+2)})$ of $C$,
where $\tilde{b}$ is the $k=1$ part of $b$.
Hence we have
$$ \rho_{w-r+2}^{(r)} \circ \tau_{w-r}^{(r+1)}(b')
     = - \tau_{w-r}^{(r+2)} \circ \rho_{w-r}^{(r+1)}(b')
     = \tau_{w-r}^{(r+2)} \circ \tau_{w-r-2}^{(r+3)}(\tilde{b})
     = 0, $$
since
$\rho_{w-r+2}^{(r)} \circ \tau_{w-r}^{(r+1)}
   = - \tau_{w-r}^{(r+2)} \circ \rho_{w-r}^{(r+1)}$ and
$\tau_{w-r}^{(r+2)} \circ \tau_{w-r-2}^{(r+3)} = 0$.
By the second assertion of Lemma \ref{WeightSpectralSequenceLemma},
there exists $a \in H^{w-r}_{\text{\rm \'et}}(X^{(r)}) = E$
such that
$$ \tau_{w-r}^{(r+1)}(b')
     = \tau_{w-r}^{(r+1)} \circ \rho_{w-r}^{(r)}(a). $$
Consider the element $(b, b'-g(a)) \in B \oplus B'$
which is in the same class as $(b,b')$
in the cohomology of the second row.
Then, we have
$$ f'(b'-g(a)) = \tau_{w-r}^{(r+1)}(b') -
       \tau_{w-r}^{(r+1)} \circ \rho_{w-r}^{(r)}(a)
   = 0, $$
hence $(b, b'-g(a)) \in \Ker(e+0,\,f+f)$.
This proves the surjectivity.
\end{proof}

Therefore, it is enough to prove
Lemma \ref{WeightSpectralSequenceLemma}.
Firstly, recall the primitive decomposition of
$H^{i}_{\text{\rm \'et}}(X^{(k)})$:
$$ H^{i}_{\text{\rm \'et}}(X^{(k)})
     = \bigoplus_{j \geq 0} L^j P^{i-2j}_{\text{\rm \'et}}(X^{(k)}), $$
where $P^{i-2j}_{\text{\rm \'et}}(X^{(k)})$
is the primitive part of $H^{i-2j}_{\text{\rm \'et}}(X^{(k)})$
as in \S \ref{SectionHardLefschetz}.
Recall that $\dim X^{(k)} = n-k+1$.
The idea of our proof of Lemma \ref{WeightSpectralSequenceLemma}
is to analyze the primitive decomposition carefully.
Since $\rho, \tau$ do not preserve the primitive decomposition,
for $0 \leq i \leq \dim X^{(k+1)} = n-k$, we put
$$ P^i\,\Im{\!}^0 \rho_i^{(k)}
   := \big( \Im \rho_i^{(k)} \big) \cap P^i_{\text{\rm \'et}}(X^{(k+1)}), $$
and $P^i\,\Im{\!}^0 \rho_i^{(k)} = 0$ otherwise.
Then, we define
$$ \Im{\!}^0 \rho_i^{(k)}
      := \bigoplus_{j \geq 0} L^j P^{i-2j}\,\Im{\!}^0 \rho_{i-2j},
   \quad
   \Im{\!}^1 \rho_i^{(k)}
      := \Im \rho_i^{(k)} / \Im{\!}^0 \rho_i^{(k)}. $$
Roughly speaking, $\Im{\!}^0 \rho_i^{(k)}$ is a part
which has a decomposition according to the primitive decomposition of
$H^i_{\text{\rm \'et}}(X^{(k+1)})$, and
$\Im{\!}^1 \rho_i^{(k)}$ is a part which does not have
any decomposition compatible with the primitive decomposition of
$H^i_{\text{\rm \'et}}(X^{(k+1)})$.
Similarly for $\Im \tau_i^{(k+1)}$,
for $0 \leq i \leq \dim X^{(k)} - 2 = n-k-1$,
we put
$$ P^{i+2}\,\Im{\!}^0 \tau_i^{(k+1)}
     := \Im \tau_i^{(k+1)} \cap P_{\text{\rm \'et}}^{i+2}(X^{(k)}) $$
and $P^{i+2}\,\Im{\!}^0 \tau_i^{(k+1)} = 0$ otherwise.
Then, we define
$$
 \Im{\!}^0 \tau_i^{(k+1)}
     := \bigoplus_{j \geq 0} L^j P^{i+2-2j}\,\Im{\!}^0 \tau_{i-2j}^{(k+1)},
   \quad
   \Im{\!}^1 \tau_i^{(k+1)}
     := \Im \tau_i^{(k+1)} / \Im{\!}^0 \tau_i^{(k+1)}.
$$

\begin{lemma}[Hard Lefschetz theorem for $\Im{\!}^0$]
\label{HardLefschetzIm0}
We have the primitive decomposition for $\Im{\!}^0$
in the following sense:
\begin{align*}
\Im{\!}^0 \rho_i^{(k)} =& \bigoplus_{j \geq 0}
 \big( L^j P^{i-2j}_{\text{\rm \'et}}(X^{(k+1)}) \cap
   \Im{\!}^0 \rho_i^{(k)} \big), \\
\Im{\!}^0 \tau_i^{(k+1)} &= \bigoplus_{j \geq 0}
 \big( L^j P^{i-2j}_{\text{\rm \'et}}(X^{(k)}) \cap
   \Im{\!}^0 \tau_i^{(k+1)} \big).
\end{align*}
Hence we have an analogue of the hard Lefschetz theorem
for $\Im{\!}^0$ in the following form:
$$
\begin{CD}
L^{\dim X^{(k+1)} - i} = L^{(n-k) - i}
  \colon \Im{\!}^0 \rho_i^{(k)}
  @>{\cong}>> \Im{\!}^0 \rho_{2(n-k)-i}^{(k)}, \\
L^{\dim X^{(k)}-(i+2)} = L^{(n-k+1)-(i+2)}
  \colon \Im{\!}^0 \tau_i^{(k+1)}
  @>{\cong}>> \Im{\!}^0 \tau_{2(n-k+1)-i-4}^{(k+1)}.
\end{CD}
$$
\end{lemma}

\begin{proof}[Proof of Lemma \ref{HardLefschetzIm0}]
These immediately follows from the definition of
$\Im{\!}^0 \rho_i^{(k)}$,\ $\Im{\!}^0 \tau_i^{(k+1)}$.
\end{proof}

\begin{lemma}[Hard Lefschetz theorem for $\Im{\!}^1$]
\label{HardLefschetzIm1}
We also have an analogue of the hard Lefschetz theorem
for $\Im{\!}^1$ in the following sense:
$$
\begin{CD}
L^{\dim X^{(k+1)}+1-i} = L^{(n-k)+1-i} \colon
    \Im{\!}^1 \rho_i^{(k)} @>{\cong}>>
    \Im{\!}^1 \rho_{2((n-k)+1)-i}^{(k)}, \\
L^{\dim X^{(k)}+1-(i+2)} = L^{(n-k+1)+1-(i+2)} \colon
    \Im{\!}^1 \tau_i^{(k+1)} @>{\cong}>>
    \Im{\!}^1 \tau_{2((n-k+1)+1)-(i+4)}^{(k+1)}.
\end{CD}
$$
\end{lemma}

\begin{proof}[Proof of Lemma \ref{HardLefschetzIm1}]
The surjectivity follows from the hard Lefschetz theorem.
The injectivity follows from the surjectivity of the maps
$$
\begin{CD}
  L^{(n-k)+1-i} \colon \Im{\!}^0 \rho_i^{(k)} @>>>
     \Im{\!}^0 \rho_{2((n-k)+1)-i}^{(k)} \\
  L^{(n-k+1)+1-(i+2)} \colon \Im{\!}^0 \tau_i^{(k+1)} @>>>
     \Im{\!}^0 \tau_{2((n-k+1)+1)-(i+4)}^{(k+1)},
\end{CD}
$$
which follows from Lemma \ref{HardLefschetzIm0}.
\end{proof}

\begin{remark}
The primitive decomposition of
$\Im{\!}^0 \rho_i^{(k)}$ (resp. $\Im{\!}^0 \tau_i^{(k+1)}$)
is the same as that of $H^{i}_{\text{\rm \'et}}(X^{(k+1)})$
(resp. $H^{i+2}_{\text{\rm \'et}}(X^{(k)})$).
On the other hand, the primitive decomposition of
$\Im{\!}^1 \rho_i^{(k)}$ (resp. $\Im{\!}^1 \tau_i^{(k+1)}$)
is the same as that of $H^{i}_{\text{\rm \'et}}(X^{(k)})$
(resp. $H^{i}_{\text{\rm \'et}}(X^{(k+1)})$).
Hence the center of symmetry of $\Im{\!}^1$
differs by 1 from that of $\Im{\!}^0$.
\end{remark}

\begin{lemma}
\label{LemmaDuality}
For all $i$, we have the following equality of dimensions:
$$
  \dim_{\Q_l} \Im{\!}^0 \rho_i^{(k)}
      = \dim_{\Q_l} \Im{\!}^1 \tau_i^{(k+1)}, \qquad
  \dim_{\Q_l} \Im{\!}^1 \rho_{i+2}^{(k)}
      = \dim_{\Q_l} \Im{\!}^0 \tau_i^{(k+1)}. $$
\end{lemma}

\begin{proof}[Proof of Lemma \ref{LemmaDuality}]
Since $\rho_i^{(k)}$ and $\tau_{2(n-k)-i}^{(k+1)}$
are dual to each other up to sign (see (\ref{CupProductDuality})),
we have
$\dim_{\Q_l} \Im \rho_i^{(k)}
    = \dim_{\Q_l} \Im \tau_{2(n-k)-i}^{(k+1)}$.
Hence we have
\begin{equation}
\label{DualityEquality1}
   \dim_{\Q_l} \Im{\!}^0 \rho_i^{(k)}
        + \dim_{\Q_l} \Im{\!}^1 \rho_i^{(k)}
   = \dim_{\Q_l} \Im{\!}^0 \tau_{2(n-k)-i}^{(k+1)}
        + \dim_{\Q_l} \Im{\!}^1 \tau_{2(n-k)-i}^{(k+1)}.
\end{equation}
Since
$\dim_{\Q_l} \Im{\!}^0 \tau_{2(n-k)-i}^{(k+1)}
    = \dim_{\Q_l} \Im{\!}^0 \tau_{i-2}^{(k+1)}$
by Lemma \ref{HardLefschetzIm0},
and
$\dim_{\Q_l} \Im{\!}^1 \tau_{2(n-k)-i}^{(k+1)}
    = \dim_{\Q_l} \Im{\!}^1 \tau_{i}^{(k+1)}$
by Lemma \ref{HardLefschetzIm1},
we have
$$ \dim_{\Q_l} \Im{\!}^0 \rho_i^{(k)}
        + \dim_{\Q_l} \Im{\!}^1 \rho_i^{(k)}
   = \dim_{\Q_l} \Im{\!}^0 \tau_{i-2}^{(k+1)}
        + \dim_{\Q_l} \Im{\!}^1 \tau_i^{(k+1)}. $$
Equivalently, we have
\begin{equation}
\label{DualityEquality2}
   \dim_{\Q_l} \Im{\!}^0 \rho_i^{(k)}
        - \dim_{\Q_l} \Im{\!}^1 \tau_i^{(k+1)}
   = - \big( \dim_{\Q_l} \Im{\!}^1 \rho_i^{(k)}
        - \dim_{\Q_l} \Im{\!}^0 \tau_{i-2}^{(k+1)} \big).
\end{equation}

On the other hand, we have
$\dim_{\Q_l} \Im{\!}^0 \rho_{i}^{(k)}
    = \dim_{\Q_l} \Im{\!}^0 \rho_{2(n-k)-i}^{(k)}$
by Lemma \ref{HardLefschetzIm0},
and
$\dim_{\Q_l} \Im{\!}^1 \rho_{i}^{(k)}
    = \dim_{\Q_l} \Im{\!}^1 \rho_{2(n-k)-i+2}^{(k)}$
by Lemma \ref{HardLefschetzIm1}.
Hence, by (\ref{DualityEquality1}), we have
\begin{align*}
  & \dim_{\Q_l} \Im{\!}^0 \rho_{2(n-k)-i}^{(k)}
        + \dim_{\Q_l} \Im{\!}^1 \rho_{2(n-k)-i+2}^{(k)} \\
  &\hspace*{1.2in}= \dim_{\Q_l} \Im{\!}^0 \tau_{2(n-k)-i}^{(k+1)}
        + \dim_{\Q_l} \Im{\!}^1 \tau_{2(n-k)-i}^{(k+1)}.
\end{align*}
From this, by replacing $2(n-k)-i$ by $i$, we have
\begin{equation}
\label{DualityEquality3}
  \dim_{\Q_l} \Im{\!}^0 \rho_{i}^{(k)}
     - \dim_{\Q_l} \Im{\!}^1 \tau_{i}^{(k+1)}
   = - \big( \dim_{\Q_l} \Im{\!}^1 \rho_{i+2}^{(k)}
     - \dim_{\Q_l} \Im{\!}^0 \tau_{i}^{(k+1)} \big).
\end{equation}

By combining (\ref{DualityEquality2}) and
(\ref{DualityEquality3}), we have
$$ \dim_{\Q_l} \Im{\!}^1 \rho_i^{(k)}
     - \dim_{\Q_l} \Im{\!}^0 \tau_{i-2}^{(k+1)}
   = \dim_{\Q_l} \Im{\!}^1 \rho_{i+2}^{(k)}
     - \dim_{\Q_l} \Im{\!}^0 \tau_{i}^{(k+1)}. $$
Hence the difference
$\dim_{\Q_l} \Im{\!}^1 \rho_{i+2}^{(k)}
     - \dim_{\Q_l} \Im{\!}^0 \tau_i^{(k+1)}$
is independent of $i$.
Since this is zero for sufficiently large (or small) $i$,
we have
$$ \dim_{\Q_l} \Im{\!}^1 \rho_{i+2}^{(k)}
     - \dim_{\Q_l} \Im{\!}^0 \tau_i^{(k+1)} = 0
   \qquad \text{for all $i$}, $$
which proves the second assertion.
Hence we also have the first assertion
by (\ref{DualityEquality2}) or (\ref{DualityEquality3}).
\end{proof}

We call the pairing on $H^{i}_{\text{\rm \'et}}(X^{(k)})$
defined by
$$ H^{i}_{\text{\rm \'et}}(X^{(k)}) \times H^{i}_{\text{\rm \'et}}(X^{(k)})
       \to \Q_l, \qquad
   (a,b) \mapsto L^{n-k+1-i} a \cup b $$
the {\em sum of the Lefschetz pairings}.
We denote it simply by $\langle,\rangle$.
The pairing $\langle,\rangle$ is nondegenerate
by the hard Lefschetz theorem (\S \ref{SectionHardLefschetz}).
However, the following lemma is far from trivial,
and we use the assumptions of Proposition \ref{WMCHodgeProp}
to prove it.

\begin{lemma}
\label{LemmaNondegeneracy}
The restriction of $\langle,\rangle$
on $H^i_{\text{\rm \'et}}(X^{(k+1)})$ to $\Im{\!}^0 \rho_i^{(k)}$
is nondegenerate. Similarly,
The restriction of $\langle,\rangle$
on $H^{i+2}_{\text{\rm \'et}}(X^{(k)})$ to $\Im{\!}^0 \tau_i^{(k+1)}$
is also nondegenerate.
\end{lemma}

\begin{proof}[Proof of Lemma \ref{LemmaNondegeneracy}]
Since the decomposition of $\Im{\!}^0 \rho_i^{(k)}$
in Lemma \ref{HardLefschetzIm0} is orthogonal with respect
to $\langle,\rangle$,
it is enough to show that
the restriction of $\langle,\rangle$
to $L^j P^{i-2j}\,\Im{\!}^0 \rho_i^{(k)}$
is nondegenerate for each $j$.
By the assumptions of Proposition \ref{WMCHodgeProp},
we observe that all cohomology groups of $X^{(k)}, X^{(k+1)}$,
maps $\rho, \tau, L$, and pairings $\langle,\rangle$
have natural $\Q$-structures (see \S \ref{SectionNotation}).
In particular, the inclusion
$$ L^j P^{i-2j}\,\Im{\!}^0 \rho_i^{(k)}
       \subset L^j P^{i-2j}_{\text{\rm \'et}}(X^{(k+1)}) $$
has a $\Q$-structure.
Namely, let $V \subset L^j P^{i-2j}\,\Im{\!}^0 \rho_i^{(k)}$ and
$W \subset L^j P^{i-2j}_{\text{\rm \'et}}(X^{(k+1)})$
be the $\Q$-vector subspaces generated by algebraic cycles,
then we have
$$ V \otimes_{\Q} \Q_l = L^j P^{i-2j}\,\Im{\!}^0 \rho_i^{(k)},
   \qquad
   W \otimes_{\Q} \Q_l = L^j P^{i-2j}_{\text{\rm \'et}}(X^{(k+1)}). $$
It is enough to show that
the restriction of $\langle,\rangle$ to $V$ is nondegenerate.
We may assume $i$ is even
because there is no odd degree cohomology.
By assumption, the Hodge standard conjecture
holds for all irreducible components of $X^{(k+1)}$.
Hence the restriction of $\langle,\rangle$ to $W$ is positive
(resp. negative) definite if $(i-2j)/2$ is even (resp. odd).
Therefore, the restriction of $\langle,\rangle$ to $V$
is nondegenerate because it is definite.
\end{proof}

\begin{lemma}
\label{LemmaIsomorphism}
The composition of the following maps is an isomorphism:
$$
\begin{CD}
\Im{\!}^0 \rho_i^{(k)} @>{\tau_i^{(k+1)}}>> \Im \tau_i^{(k+1)}
  @>>> \Im{\!}^1 \tau_i^{(k+1)} =
           \Im \tau_i^{(k+1)} / \Im{\!}^0 \tau_i^{(k+1)}.
\end{CD}
$$
Similarly, the composition of the following maps is also
an isomorphism:
$$
\begin{CD}
\Im{\!}^0 \tau_i^{(k+1)} @>{\rho_{i+2}^{(k)}}>> \Im \rho_{i+2}^{(k)}
  @>>> \Im{\!}^1 \rho_{i+2}^{(k)} =
           \Im \rho_{i+2}^{(k)} / \Im{\!}^0 \rho_{i+2}^{(k)}.
\end{CD}
$$
\end{lemma}

\begin{proof}[Proof of Lemma \ref{LemmaIsomorphism}]
By Lemma \ref{LemmaDuality},
we have only to show that the compositions are injective.
Assume that the composition of the first row is not injective
and take a nonzero $x \in \Im{\!}^0 \rho_{i}^{(k)}$
such that $\tau_i^{(k+1)}(x) \in \Im{\!}^0 \tau_{i}^{(k+1)}$.
If $\tau_i^{(k+1)}(x) \neq 0$, by Lemma \ref{LemmaNondegeneracy},
there exists $\tau_{i}^{(k+1)}(y) \in \Im{\!}^0 \tau_i^{(k+1)}$
such that
$$ L^{2(n-k+1)-(i+2)}\,\tau_i^{(k+1)}(x) \cup \tau_{i}^{(k+1)}(y) \neq 0. $$
However, since $x = \rho_{i}^{(k)}(x')$ for some
$x' \in H^{i}_{\text{\rm \'et}}(X^{(k)})$, by (\ref{CupProductDuality}),
we have
\begin{align*}
   0 &\neq L^{(n-k+1)-(i+2)}\,\tau_i^{(k+1)}(x)
               \cup \tau_{i}^{(k+1)}(y) \\
     &= \pm L^{(n-k+1)-(i+2)} \rho_{i+2}^{(k)} \circ \tau_i^{(k+1)}
               \circ \rho_{i}^{(k)}(x') \cup y \\
     &= 0.
\end{align*}
Recall that $\rho_{i+2}^{(k)} \circ \tau_i^{(k+1)} \circ \rho_{i}^{(k)} = 0$.
This is a contradiction.
Hence we have $\tau_i^{(k+1)}(x) = 0$.
By Lemma \ref{LemmaNondegeneracy} again,
there exists $\rho_{i}^{(k)}(z) \in \Im{\!}^0 \rho_{i}^{(k)}$
such that
$$ L^{n-k-i}\,x \cup \rho_{i}^{(k)}(z) \neq 0. $$
However, by (\ref{CupProductDuality}), we have
$$ 0 \neq L^{n-k-i}\,x \cup \rho_{i}^{(k)}(z)
     =  \pm L^{n-k-i}\,\tau_{i}^{(k+1)}(x) \cup z
     = 0 $$
since $\tau_i^{(k+1)}(x) = 0$.
This is a contradiction.
Hence the composition of the first row is injective.
Similarly, we see that
the composition of the second row is also injective.
\end{proof}

\begin{lemma}
\label{LemmaOrthogonal}
By Lemma \ref{LemmaIsomorphism},
the surjections
$\Im \tau_i^{(k+1)} \to \Im{\!}^1 \tau_i^{(k+1)}$ and
$\Im \rho_{i+2}^{(k)} \to \Im{\!}^1 \rho_{i+2}^{(k)}$
have canonical splittings.
Therefore, we have the following decompositions of
$\Im \tau_i^{(k+1)}$ and $\Im \rho_{i+2}^{(k)}$:
\begin{equation}
\label{DecompositionIm}
   \Im \tau_i^{(k+1)} = \Im{\!}^0 \tau_i^{(k+1)}
         \oplus \Im{\!}^1 \tau_i^{(k+1)}, \quad
   \Im \rho_{i+2}^{(k)} = \Im{\!}^0 \rho_{i+2}^{(k)}
         \oplus \Im{\!}^1 \rho_{i+2}^{(k)}.
\end{equation}
These are orthogonal decompositions
with respect to $\langle,\rangle$.
\end{lemma}

\begin{proof}[Proof of Lemma \ref{LemmaOrthogonal}]
The proof is easy and immediate.
Take
$(a,b) \in \Im{\!}^0 \tau_i^{(k+1)} \oplus \Im{\!}^1 \tau_i^{(k+1)}$.
Then, by the first assertion of Lemma \ref{LemmaIsomorphism},
we can write
$a = \tau_i^{(k+1)}(a')$ and $b = \tau_i^{(k+1)} \circ \rho_i^{(k)}(b')$
for some $a', b'$.
Then, by (\ref{CupProductDuality}), we have
\begin{align*}
   L^{n-k+1-(i+2)}\,a \cup b
     &= L^{n-k+1-(i+2)}\,\tau_i^{(k+1)}(a')
         \cup \big( \tau_i^{(k+1)} \circ \rho_i^{(k)}(b')\big) \\
     &= \pm L^{n-k+1-(i+2)}\,a' \cup \big( \rho_{i+2}^{(k+1)}
         \circ \tau_i^{(k+1)} \circ \rho_i^{(k)}(b') \big) \\
     &= 0,
\end{align*}
since $\rho_{i+2}^{(k+1)} \circ \tau_i^{(k+1)}
  \circ \rho_i^{(k)} = 0$.
Similarly, we see that
the decomposition
$\Im \rho_{i+2}^{(k)} = \Im{\!}^0 \rho_{i+2}^{(k)}
     \oplus \Im{\!}^1 \rho_{i+2}^{(k)}$
is also orthogonal with respect to $\langle,\rangle$.
\end{proof}

Finally, by combining above results,
we shall prove Lemma \ref{WeightSpectralSequenceLemma}.

\begin{proof}[Proof of Lemma \ref{WeightSpectralSequenceLemma}]
Firstly, we shall prove the first assertion.
Since
$\tau_i^{(k+1)} \circ \rho_i^{(k)} \circ \tau_{i-2}^{(k+1)} = 0$,
the inclusion
$$ \big( \Ker \tau_i^{(k+1)} \big) \cap \big( \Im \rho_i^{(k)} \big)
   \supset \Im \big( \rho_i^{(k)} \circ \tau_{i-2}^{(k+1)} \big) $$
is obvious.
To prove the opposite inclusion,
we take $x \in \big( \Ker \tau_i^{(k+1)} \big)
   \cap \big( \Im \rho_i^{(k)} \big)$.
Then, for all $y \in H^i_{\text{\rm \'et}}(X^{(k)})$,
by (\ref{CupProductDuality}), we have
$$
   L^{n-k-i}\,x \cup \rho_i^{(k)}(y)
      = \pm L^{n-k-i}\,\tau_{i}^{(k+1)}(x) \cup y
      = 0
$$
since $\tau_{i}^{(k+1)}(x) = 0$.
By Lemma \ref{LemmaOrthogonal}, we write $x = x_0 + x_1$
for $x_0 \in \Im{\!}^0 \rho_i^{(k)}$ and
$x_1 \in \Im{\!}^1 \rho_i^{(k)}$.
If $x_0 \neq 0$, by Lemma \ref{LemmaNondegeneracy}
and Lemma \ref{LemmaOrthogonal},
there exists $y' \in \Im{\!}^0 \rho_i^{(k)}$
such that
\begin{align*}
   L^{n-k-i}\,x \cup y'
     &= L^{n-k-i}\,x_0 \cup y' + L^{n-k-i}\,x_1 \cup y' \\
     &= L^{n-k-i}\,x_0 \cup y' \\
     &\neq 0,
\end{align*}
which is absurd. Hence $x_0 = 0$.
Then, we have
$$ x = x_1 \in \Im{\!}^1 \rho_i^{(k)}
   = \rho_i^{(k)}(\Im{\!}^0 \tau_{i-2}^{(k+1)})
   \subset \Im \big( \rho_i^{(k)} \circ \tau_{i-2}^{(k+1)} \big). $$
This proves the first assertion.
Similarly, we can prove the second assertion
by using Lemma \ref{LemmaNondegeneracy}
and Lemma \ref{LemmaOrthogonal}.
Therefore, Lemma \ref{WeightSpectralSequenceLemma} is proved,
and the proof of Proposition \ref{WMCHodgeProp} is complete.
\end{proof}

\begin{remark}
In the above proof,
Lemma \ref{LemmaNondegeneracy}
is the only part where we crucially use
the assumptions of Proposition \ref{WMCHodgeProp}.
In \cite{MSaito}, M. Saito proved a corresponding statement
by using polarized Hodge structures.
\end{remark}

\section{Application to $p$-adically uniformized varieties}
\label{SectionP-adic}

In this section, we give a proof of the main theorem
(Theorem \ref{MainTheorem}) of this paper.
We also give some immediate applications.

\subsection{Drinfeld upper half spaces and
$p$-adic uniformization}
\label{DrinfeldUpperHalfSpace}

Let $K$ be a complete discrete valuation field
with finite residue field $\F_q$ of characteristic $p>0$.
Let $\widehat{\Omega}^d_{K}$ be the Drinfeld upper half space
of dimension $d \geq 1$ over $K$, which is a rigid analytic space
obtained by removing all $K$-rational hyperplanes from $\P^{d}_K$.
We have a natural action of $\PGL_{d+1}(K)$ on $\widehat{\Omega}^d_{K}$.

As a formal scheme, $\widehat{\Omega}^d_{K}$ can be constructed
as follows. Take the projective space $\P^d_{\O_K}$ over $\O_K$.
Then, take successive blowing-ups of $\P^d_{\O_K}$ along
all linear subvarieties in the special fiber $\P^d_{\F_q}$
as in \S \ref{SubsectionConstruction}.
By continuing this process for {\em all} exceptional divisors
appearing in the blowing-ups,
we obtain a formal scheme $\widehat{\Omega}^d_{\O_K}$
locally of finite type over $\text{Spf}\,\O_K$. 
By construction, the rigid analytic space
associated with $\widehat{\Omega}^d_{\O_K}$ is
isomorphic to $\widehat{\Omega}^d_{K}$.
In other words, $\widehat{\Omega}^d_{K}$ is
the \lq\lq generic fiber" in the sense of Raynaud
(\cite{Raynaud}) of the formal scheme
$\widehat{\Omega}^d_{\O_K}$ over $\text{Spf}\,\O_K$.

Let $\Gamma$ be a cocompact torsion free discrete subgroup
of $\PGL_{d+1}(K)$.
We have a natural action of $\Gamma$
on $\widehat{\Omega}^d_{\O_K}$,
and this action is discontinuous with respect to Zariski topology.
Then, we take a quotient
$ \widehat{\X}_{\Gamma} :=
  \Gamma \backslash \widehat{\Omega}^d_{\O_K} $
as a formal scheme.
Mustafin, Kurihara showed that the relative dualizing
sheaf $\omega_{\widehat{\X}_{\Gamma}/\O_K}$ is
invertible and ample.
Therefore, $\widehat{\X}_{\Gamma}$ can
be algebraized to a projective scheme
$\X_{\Gamma}$ over $\O_K$.
The generic fiber
$X_{\Gamma} := \X_{\Gamma} \otimes_{\O_K} K$
is a projective smooth variety over $K$
whose associated rigid analytic space is
the rigid analytic quotient
$\Gamma \backslash \widehat{\Omega}^d_{K}$.
By construction, $\X_{\Gamma}$ is
a proper semistable model of $X_{\Gamma}$ over $\O_K$
whose special fiber is described by the cell complex
$\Gamma \backslash {\mathfrak T}$,
where $\mathfrak T$ denotes the Bruhat-Tits building of
$\PGL_{d+1}(K)$
(for details, see \cite{Mustafin}, Theorem 4.1,
\cite{Kurihara}, Theorem 2.2.5).

\subsection{Proof of the main theorem (Theorem \ref{MainTheorem})}
\label{SectionProofMainTheorem}

As in \S \ref{DrinfeldUpperHalfSpace},
let $K$ be a complete discrete valuation field
with finite residue field $\F_q$ of characteristic $p>0$,
$\Gamma \subset \PGL_{d+1}(K)$
a cocompact torsion free discrete subgroup,
and $\X_{\Gamma}$ an algebraization of
$\Gamma \backslash \widehat{\Omega}^d_{\O_K}$
with generic fiber $X_{\Gamma}$
as in \S \ref{DrinfeldUpperHalfSpace}.

The following lemma seems well-known.

\begin{lemma}
\label{LemmaEtaleCover}
Let $X, Y$ be proper smooth variety of dimension $d$
over $K$, and $f \colon X \to Y$ a proper surjective
generically finite morphism of degree $\text{\rm deg}\,f$
($f$ is called an {\rm alteration} in \cite{deJong}).
If the weight-monodromy conjecture
(Conjecture \ref{WMC}) holds for $X$,
then Conjecture \ref{WMC} holds also for $Y$.
\end{lemma}

\begin{proof}
Since $f$ is finite of degree $\text{\rm deg}\,f$
on some open dense subscheme $U \subset X$,
$f^{\ast} \colon H^{k}_{\text{\rm \'et}}(Y_{\overline{K}}, \Q_l)
   \to H^{k}_{\text{\rm \'et}}(X_{\overline{K}}, \Q_l)$
and
$f_{\ast} \colon H^{k}_{\text{\rm \'et}}(X_{\overline{K}}, \Q_l)
   \to H^{k}_{\text{\rm \'et}}(Y_{\overline{K}}, \Q_l)$
satisfy $f_{\ast} \circ f^{\ast} =$
multiplication by $\text{deg}\, f$.
Since $f_{\ast},\ f^{\ast}$ commute with
$\Gal(\overline{K}/K)$-action,
$H^{k}_{\text{\rm \'et}}(Y_{\overline{K}}, \Q_l)$
is a direct summand of
$H^{k}_{\text{\rm \'et}}(X_{\overline{K}}, \Q_l)$
as a $\Gal(\overline{K}/K)$-representation.
Therefore, the restriction of the monodromy filtration
(resp. weight filtration)
on $H^{k}_{\text{\rm \'et}}(X_{\overline{K}}, \Q_l)$
to $H^{k}_{\text{\rm \'et}}(Y_{\overline{K}}, \Q_l)$
coincides with the monodromy filtration
(resp. weight filtration) on
$H^{k}_{\text{\rm \'et}}(Y_{\overline{K}}, \Q_l)$ itself.
Hence the assertion follows.
\end{proof}

Now, we shall prove the main theorem (Theorem \ref{MainTheorem})
of this paper.

\begin{proof}[Proof of Theorem \ref{MainTheorem}]
First of all, it is known that the intersection of
all finite index subgroups
$\Gamma' \subset \Gamma$ is equal to the identity
(\cite{Garland}, Theorem 2.7, Proof, (b)).
From this, by looking at the action of $\Gamma$ on
the Bruhat-Tits building $\mathfrak T$ of $\PGL_{d+1}(K)$
(\cite{Mustafin}, \cite{Kurihara}),
it is easy to see that
$\X_{\Gamma'}$ is {\em strictly} semistable
for some finite index subgroup $\Gamma' \subset \Gamma$.
Since we have a finite \'etale covering
$X_{\Gamma'} \to X_{\Gamma}$, by Lemma \ref{LemmaEtaleCover},
we may assume that $\X_{\Gamma}$ is a strictly semistable model
of $X_{\Gamma}$ over $\O_K$.

Under this assumption,
from the construction of $\X_{\Gamma}$
in \S \ref{DrinfeldUpperHalfSpace},
we see that all irreducible components
$X_1,\ldots,X_m$ of the special fiber of
$\X_{\Gamma}$ are isomorphic to the variety $B^d$
constructed in \S \ref{SubsectionConstruction}.
By Proposition \ref{PropCombinatorics}, 4, 5,
for $i \neq j$, each irreducible component of
$X_i \cap X_j$ is isomorphic to a divisor
of the form $D_{V}$ on $B^d$ in \S \ref{SubsectionDivisors}.
Moreover, by induction,
we see that, for $1 \leq i_1 < \cdots < i_k \leq m$,
each irreducible component $Y$ of
$X_{i_1} \cap \cdots \cap X_{i_k}$
is isomorphic to the product
$$ Y \cong B^{n_1} \times \cdots \times B^{n_k} $$
for some $n_1,\ldots,n_k \geq 1$ with
$n_1 + \cdots + n_k = d-k+1$.

Let $\L$ be the relative dualizing sheaf
$\omega_{\X_{\Gamma}/\O_K}$, which is
invertible and ample by Mustafin, Kurihara
(\cite{Mustafin}, \cite{Kurihara}).
Take $i$ with $1 \leq i \leq m$ and
fix an isomorphism $X_i \cong B^d$.
Then, by an explicit calculation in
\cite{Mustafin}, \cite{Kurihara},
the restriction of $\L$ to each $X_i$ is isomorphic to
    $$ -(n+1) f^{\ast} \O_{\P^d}(1) + \sum_{d=0}^{n-1} (n-d) D_d $$
(for notation, see \S \ref{SubsectionConstruction},
 \S \ref{SubsectionDivisors},
 see also Remark \ref{RemarkAmpleDivisorP-adicUniformization}),
which is an ample $\PGL_{d+1}(\F_q)$-invariant
divisor on $X_i \cong B^d$.
Therefore, by Proposition \ref{MainPropLefschetzHodge},
the Hodge standard conjecture (Conjecture \ref{HodgeConj})
holds for $(X_i,\L|_{X_i})$.

Moreover, for $1 \leq i_1 < \cdots < i_k \leq m$,
take an irreducible component $Y$ of
$X_{i_1} \cap \cdots \cap X_{i_k}$.
Then, $Y$ is isomorphic to
$B^{n_1} \times \cdots \times B^{n_k}$
for some $n_1,\ldots,n_k \geq 1$ with
$n_1 + \cdots + n_k = d-k+1$.
By the construction of $B^n$ in
\S \ref{SubsectionConstruction} and the K\"unneth
formula, $Y$ satisfies
Assumption \ref{AssumptionCohomology}
in \S \ref{SectionNotation}
(see also \S \ref{SubsectionProductVarieties}).
Let $\pr_j \colon Y \to B^{n_j}$ be the projection
to the $i$-th factor for $j = 1,\ldots,k$.
Then,
by applying Proposition \ref{PropCombinatorics}, 4, 5
and Proposition \ref{PGLinvariantProp}, 4
inductively,
we see that the restriction of $\L$ to $Y$
is of the form
    $$ \L|_Y = \pr_1^{\ast} \O(L_1) + \cdots + \pr_k^{\ast} \O(L_k)
       \qquad \text{in} \quad H^2(Y) $$
where $L_j$ is an ample $\PGL_{n_j+1}(\F_q)$-invariant
divisor on $B^{n_j}$ for $j = 1,\ldots,k$.
By Proposition \ref{MainPropLefschetzHodge},
the Hodge standard conjecture (Conjecture \ref{HodgeConj})
holds for $(B^{n_j},L_j)$ for $j = 1,\ldots,k$.
By applying Proposition \ref{ConjectureProduct} inductively,
we conclude that
the Hodge standard conjecture (Conjecture \ref{HodgeConj})
also holds for $(Y,\L|_{Y})$.

Therefore, $\X_{\Gamma}$ satisfies all assumptions
in Proposition \ref{WMCHodgeProp}.
Hence the weight-monodromy conjecture (Conjecture \ref{WMC})
holds for $X_{\Gamma}$, and the proof of Theorem \ref{MainTheorem}
is complete.
\end{proof}

\begin{remark}
\label{MainTheoremCharP}
Theorem \ref{MainTheorem} is new only in
mixed characteristic.
However, we note that Theorem \ref{MainTheorem}
in characteristic $p>0$ does not automatically follow
from Deligne's results
in \cite{WeilII} because $\X_{\Gamma}$ does not
automatically come from a family of varieties
over a curve over a finite field
(This was explained by Illusie in \cite{Illusie3}, 8.7).
We need a specialization argument to reduce the general
characteristic $p>0$ case to Deligne's case
(\cite{Terasoma}, \cite{Ito1}).
\end{remark}

\subsection{Schneider-Stuhler's conjecture on the filtration $F^{\bullet}$}
\label{SectionSchneiderStuhlerConjecture}

In \cite{Schneider-Stuhler},
Schneider-Stuhler computed the cohomology of
$\widehat{\Omega}^d_{K}$ and its quotient
$\Gamma \backslash \widehat{\Omega}^d_{K}$
by methods from rigid analytic geometry
and representation theory.
Their results are valid for cohomology theories
satisfying certain axioms,
and they showed that de Rham cohomology
and rigid analytic $l$-adic cohomology
satisfy these axioms.
Precisely speaking, at that time,
they used unpublished results of Gabber
to check the axioms for rigid analytic $l$-adic cohomology.
After that, these axioms were established
by Berkovich (\cite{Berkovich}),
de Jong\,-\,van der Put (\cite{deJong-vanderPut}).
By comparison theorem between
rigid analytic and usual $l$-adic cohomology,
Schneider-Stuhler's computation is valid for
$l$-adic cohomology of the variety $X_{\Gamma}$.
Here we combine our results
with the results of Schneider-Stuhler in
\cite{Schneider-Stuhler},
and prove Schneider-Stuhler's conjecture
on the filtration $F^{\bullet}$ for
$l$-adic cohomology.

First of all, we recall the results of Schneider-Stuhler
(for details, see \cite{Schneider-Stuhler}, Theorem 4).
As in \S \ref{DrinfeldUpperHalfSpace},
let $K$ be a complete discrete valuation field
with finite residue field $\F_q$ of characteristic $p>0$.
For a cocompact torsion free discrete subgroup
$\Gamma \subset \PGL_{d+1}(K)$,
let
$$ \text{Ind}_{\Gamma} := C^{\infty} \big( \PGL_{d+1}(K)/\Gamma, \C \big) $$
be the $\PGL_{d+1}(K)$-representation 
induced from the trivial character on $\Gamma$.
Let $\mu(\Gamma)$ be the multiplicity of
the Steinberg representation in $\text{Ind}_{\Gamma}$.
In \cite{Schneider-Stuhler}, \S 5, Schneider-Stuhler
explicitly computed the $E_2$-terms of
a Hochschild-Serre type covering spectral sequence:
\begin{equation}
\label{CoveringSpectralSequence}
   E_2^{r,s} = H^{r} \big( \Gamma, H^s_{\text{\rm \'et}}
     (\widehat{\Omega}^d_{K} \otimes_K \overline{K},\,\Q_l) \big)
   \Longrightarrow
   H^{r+s}_{\text{\rm \'et}}(X_{\Gamma} \otimes_K \overline{K},\,\Q_l)
\end{equation}
(\cite{Schneider-Stuhler}, \S 5, Proposition 2), and proved
\begin{equation}
\label{SchneiderStuhlerComputation1}
   H^{k}_{\text{\rm \'et}}(X_{\Gamma} \otimes_K \overline{K},\,\Q_l)
     \cong
     \begin{cases}
       \Q_l\!\left( -\frac{k}{2} \right)
           & \text{if $k$ is even,\ $0 \leq k \leq 2d$,\ $k \neq d$} \\
       0   & \text{if $k$ is odd,\ $k \neq d$.}
     \end{cases}
\end{equation}
For the middle degree $k = d$, they proved that
the covering spectral sequence
(\ref{CoveringSpectralSequence}) defines
a decreasing filtration $F^{\bullet}$ on
$V := H^{d}_{\text{\rm \'et}}(X_{\Gamma} \otimes_K \overline{K},\,\Q_l)$
such that
$$  V = F^0 V \supset F^1 V \supset \cdots \supset
        F^{d} V \supset F^{d+1} V = 0, $$
\begin{align}
\label{SchneiderStuhlerComputation2}
  F^r V / F^{r+1} V
  &\cong
  \begin{cases}
     \Q_l(r-d)^{\oplus \mu(\Gamma)}
         & \text{if}\ 0 \leq r \leq d,\ r \neq \frac{d}{2} \\
     \Q_l\!\left( -\frac{d}{2} \right)^{\oplus (\mu(\Gamma)+1)}
         & \text{if}\ r = \frac{d}{2} \\
     0 & \text{otherwise},
  \end{cases}
\end{align}
and conjectured that $F^{\bullet}$ essentially coincides with
the monodromy filtration (see \cite{Schneider-Stuhler},
introduction and a remark following Theorem 5).

\begin{theorem}[Schneider-Stuhler's conjecture for $l$-adic cohomology]
\label{SchneiderStuhlerConjecture}
Let $M_{\bullet}$ be the monodromy filtration on $V$
(Definition \ref{DefMonodromyFiltration}).
Define an increasing filtration $F'_{\bullet}$ on $V$ by
$ F'_{i} V = F^{- \lfloor i/2 \rfloor} V. $
Then, we have
$$ M_i\,V = F'_{i-d}\,V \qquad \text{for all}\ i. $$
\end{theorem}

\begin{proof}
We compute the graded quotients of $F'_{\bullet}$
explicitly.
If $i$ is even, we have
\begin{align*}
  F'_{i}\,V / F'_{i-1}\,V
   &= F^{- i/2} V
       / F^{- i/2 + 1} V \\
   &= \begin{cases}
       \Q_l\!\left(- \frac{i}{2} - d \right)^{\oplus \mu(\Gamma)}
           & \text{if}\ -2d \leq i \leq 0,\ i \neq -d \\
       \Q_l\!\left( -\frac{d}{2} \right)^{\oplus (\mu(\Gamma)+1)}
           & \text{if}\ i = - d \\
       0 & \text{otherwise}.
     \end{cases}
\end{align*}
If $i$ is odd, we have
$\lfloor i/2 \rfloor = \lfloor (i-1)/2 \rfloor$
hence $F'_{i}\,V / F'_{i-1}\,V = 0$.
Therefore, each
$F'_{i}\,V / F'_{i-1}\,V$ has weight $i + 2d$
as a $\Gal(\overline{\F}_q/\F_q)$-representation.
Now the assertion follows from Theorem \ref{MainTheorem}.
\end{proof}

\subsection{Application to the local zeta function of $X_{\Gamma}$}
\label{SubsectionLocalZetaFunctions}

Let notation be as in \S \ref{SectionSchneiderStuhlerConjecture}.
Let us recall the definition of the local zeta functions
(for details, see \cite{Serre2}, \cite{Rapoport}).
For a continuous $\Gal(\overline{K}/K)$-representation $V$
over $\Q_l$, the {\em local $L$-function} of $V$ is defined by
  $$ L(s,V) := \det \big( 1 - q^{-s} \cdot \Fr_q\,;\,V^{I_K} \big)^{-1} $$
(for notation, see \S \ref{SectionMonodromyFiltration},
 \S \ref{SectionWeightFiltration}).
For a variety $X$ over $K$, the {\em local zeta function} of $X$
is defined by
$$ \zeta(s,X) := \prod_{k=0}^{2 \dim X}
   L \big( s,\,H^{k}_{\text{\rm \'et}}(X_{\overline{K}},\Q_l) \big)^{(-1)^k}. $$

Let us compute the local zeta function of $X_{\Gamma}$.
For an even integer $k$ with $0 \leq k \leq 2d,\ k \neq d$,
by (\ref{SchneiderStuhlerComputation1}),
we easily see that
$$ L \big( s, H^{k}_{\text{\rm \'et}}(X_{\Gamma} \otimes_K
        \overline{K},\,\Q_l) \big)
     = \det \big( 1 - q^{-s} \cdot \Fr_q\,;\,
         \Q_l{\textstyle \left( -\frac{k}{2} \right)} \big)^{-1}
     = \frac{1}{1-q^{k/2-s}}. $$
On the other hand, for the middle degree $k=d$,
we can not compute the local $L$-function only from
the results of Schneider-Stuhler (\ref{SchneiderStuhlerComputation2})
because (\ref{SchneiderStuhlerComputation2}) only gives us
the semisimplification of the $\Gal(\overline{K}/K)$-representation
$H^{d}_{\text{\rm \'et}}(X_{\Gamma} \otimes_K \overline{K},\,\Q_l)$.
We need the weight-monodromy conjecture to recover
the local $L$-function from the semisimplification
(see \cite{Rapoport}, \S 2).

Now, we shall use Theorem \ref{SchneiderStuhlerConjecture}.
By Theorem \ref{SchneiderStuhlerConjecture},
we see that the inertia fixed part is given by
$$ H^{d}_{\text{\rm \'et}}(X_{\Gamma} \otimes_K \overline{K},\,\Q_l)^{I_K} =
   \begin{cases}
     \Q_l^{\oplus \mu(\Gamma)} \oplus \Q_l\!\left( -\frac{d}{2} \right)
         & \text{if $d$ is even} \\
     \Q_l^{\oplus \mu(\Gamma)} & \text{if $d$ is odd}
   \end{cases} $$
(see also Definition \ref{DefMonodromyFiltration}, 3).
Therefore, we have
$$ L \big( s, H^{d}_{\text{\rm \'et}}(X_{\Gamma} \otimes_K
        \overline{K},\,\Q_l) \big) =
   \begin{cases}
     \displaystyle \frac{1}{(1-q^{-s})^{\mu(\Gamma)}}
           \cdot \frac{1}{1-q^{d/2-s}}
         & \text{if $d$ is even} \\
     \displaystyle \frac{1}{(1-q^{-s})^{\mu(\Gamma)}}
         & \text{if $d$ is odd.}
   \end{cases} $$

By combining above results, we have the following theorem.

\begin{theorem}[Local zeta function of $X_{\Gamma}$]
\label{LocalZetaFunction}
Let $\Gamma \subset \PGL_{d+1}(K)$ be
a cocompact torsion free discrete subgroup.
Then, the local zeta function of $X_{\Gamma}$ is given by
$$ \zeta(s,X_{\Gamma})
   = (1-q^{-s})^{\mu(\Gamma) \cdot (-1)^{d+1}}
        \cdot \prod_{k=0}^{d} \frac{1}{1-q^{k-s}}, $$
where $\mu(\Gamma)$ is the multiplicity of
the Steinberg representation in
$\text{\rm Ind}_{\Gamma}$
as in \S \ref{SectionSchneiderStuhlerConjecture}.
\end{theorem}

\begin{remark}
As a consequence, we also see that
the local zeta function $\zeta(s,X_{\Gamma})$
is independent of $l \neq p$.
In \S \ref{SubsectionP-adicAnalogue},
we consider the case $l=p$ via $p$-adic Hodge theory.
\end{remark}

\subsection{$p$-adic weight-monodromy conjecture}
\label{SubsectionP-adicAnalogue}

We recall some results in $p$-adic Hodge theory
to state a $p$-adic analogue of the weight-monodromy conjecture.

Let $K$ be a finite extension of $\Q_p$
with residue field $\F_q$,
$W(\F_q)$ the ring of Witt vectors with coefficients in $\F_q$,
$K_0$ the field of fractions of $W(\F_q)$,
and $X$ a proper smooth variety over $K$.
For simplicity, we assume that there is a proper strictly
semistable model $\X$ of $X$ over $\O_K$
(see Definition \ref{DefSemistableModel}).

For $V := H^{w}_{\text{\rm \'et}}(X_{\overline{K}}, \Q_p)$,
by the $C_{\text{\rm st}}$-conjecture proved by Tsuji,
there is a canonical isomorphism
\begin{equation}
\label{Comparisonp-adic}
  D_{\text{\rm st}}(V)
  := \big( V \otimes_{\Q_p} B_{\text{\rm st}} \big)^{\Gal(\overline{K}/K)}
  \cong  H^w_{\text{\rm log-crys}}
           \big( \X_{\F_q}^{\times}/W(\F_q)^{\times} \big)
           \otimes_{W(\F_q)} K_0
\end{equation}
where $B_{\text{\rm st}}$ is Fontaine's ring of $p$-adic periods,
$\Gal(\overline{K}/K)$ acts on $V \otimes_{\Q_p} B_{\text{\rm st}}$
diagonally,
and $H^w_{\text{\rm log-crys}} \big( \X_{\F_q}^{\times}/W(\F_q)^{\times} \big)$
is the log crystalline cohomology of the special fiber
$\X_{\F_q} = \X \otimes_{\O_K} \F_q$ endowed with
a natural log structure
(\cite{Fontaine}, \cite{HyodoKato}, \cite{Kato}, \cite{Tsuji}).
It is known that $\dim_{K_0} D_{\text{\rm st}}(V) = \dim_{\Q_p} V$.
Moreover, it is also known that the monodromy operator $N$ and
the Frobenius endomorphism $\varphi$ act on both hand sides,
and (\ref{Comparisonp-adic}) is compatible with the actions of them.
By using $N, \varphi$, we can define
the monodromy filtration $M_{\bullet}$ and
the weight filtration $W_{\bullet}$ on $D_{\text{\rm st}}(V)$
by the same way as in Definition \ref{DefMonodromyFiltration},
Definition \ref{DefWeightFiltration}.
Then, we have the following $p$-adic analogue of
the weight-monodromy conjecture.

\begin{conjecture}[$p$-adic weight-monodromy conjecture]
\label{WMCp-adic}
$$ M_i D_{\text{\rm st}}(V)
     = W_{i+w} D_{\text{\rm st}}(V) \qquad \text{for all}\ i. $$
\end{conjecture}

Let notation be as in \S \ref{SubsectionWeightSpectralSequence}.
Mokrane constructed a $p$-adic analogue of
the weight spectral sequence of Rapoport-Zink as follows:
\begin{center}
\begin{tabular}{l}
$\displaystyle E_1^{-r,\,w+r} = \bigoplus_{k \geq \max\{0,-r\}}
   H^{w-r-2k}_{\text{\rm crys}} \big( X^{(2k+r+1)}/W(\F_q) \big)(-r-k)$ \\
\hspace*{3in} $\displaystyle \Longrightarrow
\quad H^w_{\text{\rm log-crys}}
   \big( \X_{\F_q}^{\times}/W(\F_q)^{\times} \big)$,
\end{tabular}
\end{center}
where $H^{\ast}_{\text{\rm crys}}$ denotes the crystalline cohomology
(\cite{Mokrane}, \S 3.23, Th\'eor\`eme 3.32).
This spectral sequence has similar properties as the $l$-adic case
(compare with \S \ref{SubsectionWeightSpectralSequence}).
By the Weil conjecture for crystalline cohomology
(\cite{Katz-Messing}, \cite{ChiarellottoLeStum}),
this spectral sequence degenerates at $E_2$ modulo torsion,
and defines the weight filtration $W_{\bullet}$ on
$H^w_{\text{\rm log-crys}}
   \big( \X_{\F_q}^{\times}/W(\F_q)^{\times} \big)$.
Moreover, there is a monodromy operator $N$
satisfying the same properties as the $l$-adic case
($N$ coincides with $\nu$ in \cite{Mokrane}, \S 3.33).

Therefore, by (\ref{Comparisonp-adic}),
Conjecture \ref{WMCp-adic} is equivalent to the following
conjecture on the weight spectral sequence of Mokrane
(compare with Conjecture \ref{WMC_semistable}).

\begin{conjecture}[\cite{Mokrane}, Conjecture 3.27, \S 3.33]
\label{WMCp-adic_semistable}
$N^r$ induces an isomorphism
$$
\begin{CD}
N^r \colon E_2^{-r,\,w+r}(r) \otimes_{W(\F_q)} K_0
  @>{\cong}>> E_2^{r,\,w-r} \otimes_{W(\F_q)} K_0
\end{CD}
$$
for all $r,w$.
\end{conjecture}

\begin{theorem}[$p$-adic weight-monodromy conjecture for $X_{\Gamma}$]
\label{MainTheoremp-adic}
Let $\Gamma \subset \PGL_{d+1}(K)$ be
a cocompact torsion free discrete subgroup.
Then, a $p$-adic analogue of the weight-monodromy conjecture
(Conjecture \ref{WMCp-adic}) holds for $X_{\Gamma}$.
\end{theorem}

\begin{proof}
Note that $\X_{\Gamma}$ is semistable but
not necessarily strictly semistable in general.
However, as we see in \S \ref{SectionProofMainTheorem},
$\X_{\Gamma'}$ is strictly semistable for some finite index
subgroup $\Gamma' \subset \Gamma$.
Since a $p$-adic analogue of Lemma \ref{LemmaEtaleCover}
holds by the same reason,
we may assume that $\X_{\Gamma}$ is strictly semistable.
Then, we can prove Theorem \ref{MainTheoremp-adic} by the same
way as the $l$-adic case. We use the cycle map
for crystalline cohomology (\cite{Gillet-Messing}, \cite{Gros}).
\end{proof}

\begin{remark}
After this work was completed, the author was informed
that Ehud de Shalit obtained Theorem \ref{MainTheoremp-adic}
by a completely different method.
His proof relies on a combinatorial result of Alon-de Shalit
about harmonic cochains on the Bruhat-Tits buildings
(\cite{Alon-deShalit}, \cite{deShalit}).
\end{remark}

\subsection{Application to
the $p$-adic local zeta function of $X_{\Gamma}$}
\label{SubsectionP-adicAnalogueLocalZetaFunctions}

Let notation be as in \S \ref{SubsectionP-adicAnalogue}.
Recall that a finite dimensional continuous
$\Gal(\overline{K}/K)$-representation $V$ over $\Q_p$
is called {\em semistable} if
$\dim_{K_0} D_{\text{\rm st}}(V) = \dim_{\Q_p} V$
(\cite{Fontaine}).
For semistable $V$,
the {\em $p$-adic local $L$-function} of $V$ is defined by
$$ L_{\text{\rm $p$-adic}}(s,V)
      := \det \big( 1 - q^{-s} \cdot \varphi
                 \,;\,D_{\text{\rm st}}(V)^{N=0} \big)^{-1}, $$
where $D_{\text{\rm st}}(V)^{N=0}$ denotes
the kernel of $N$ acting on $D_{\text{\rm st}}(V)$.
In fact, we can define $L_{\text{\rm $p$-adic}}(s,V)$
for {\em potentially semistable} $V$,
but we omit it here (for details, see \cite{Fontaine}).
For a variety $X$ over $K$ such that
$H^{k}_{\text{\rm \'et}}(X_{\overline{K}},\Q_p)$ are semistable
for all $k$,
we define the {\em $p$-adic local zeta function}
$\zeta_{\text{\rm $p$-adic}}(s,X)$ of $X$ as follows:
$$ \zeta_{\text{\rm $p$-adic}}(s,X) := \prod_{k=0}^{2 \dim X}
   L_{\text{\rm $p$-adic}} \big( s,\,
     H^{k}_{\text{\rm \'et}}(X_{\overline{K}},\Q_p) \big)^{(-1)^i}. $$
Here we put the subscript \lq\lq $p$-adic'' to distinguish it from
the $l$-adic case.

\begin{theorem}[$p$-adic local zeta function of $X_{\Gamma}$]
\label{LocalZetaFunctionP-adic}
Let $\Gamma \subset \PGL_{d+1}(K)$ be
a cocompact torsion free discrete subgroup.
Then, for all $k$,
the $p$-adic \'etale cohomology
$H^{k}_{\text{\rm \'et}}(X_{\Gamma}
        \otimes_K \overline{K},\,\Q_p)$
is a semistable $\Gal(\overline{K}/K)$-representation,
and the $p$-adic local $L$-function
$L_{\text{\rm $p$-adic}} \big( s, H^{k}_{\text{\rm \'et}}(X_{\Gamma}
        \otimes_K \overline{K},\,\Q_p) \big)$
is the same as the usual local $L$-function
defined by $l$-adic cohomology for $l \neq p$
(see \S \ref{SubsectionLocalZetaFunctions}):
$$ L_{\text{\rm $p$-adic}} \big( s, H^{k}_{\text{\rm \'et}}(X_{\Gamma}
        \otimes_K \overline{K},\,\Q_p) \big)
  = L \big( s, H^{k}_{\text{\rm \'et}}(X_{\Gamma}
        \otimes_K \overline{K},\,\Q_l) \big). $$
In particular, we have the following equality of local zeta functions
of $X_{\Gamma}$:
$$ \zeta_{\text{\rm $p$-adic}}(s,X_{\Gamma})
   = \zeta(s,X_{\Gamma})
   = (1-q^{-s})^{\mu(\Gamma) \cdot (-1)^{d+1}}
        \cdot \prod_{k=0}^{d} \frac{1}{1-q^{k-s}}, $$
where $\mu(\Gamma)$ is the multiplicity of
the Steinberg representation in
$\text{\rm Ind}_{\Gamma}$
as in \S \ref{SectionSchneiderStuhlerConjecture}.
\end{theorem}

\begin{proof}
If $\X_{\Gamma}$ is strictly semistable,
$H^{k}_{\text{\rm \'et}}(X_{\Gamma}
        \otimes_K \overline{K},\,\Q_p)$
is semistable by 
the $C_{\text{\rm st}}$-conjecture proved by Tsuji
(\cite{Tsuji}).
The general case follows from this
because $\X_{\Gamma'}$ is strictly semistable for
some finite index subgroup $\Gamma' \subset \Gamma$
(see also Lemma \ref{LemmaEtaleCover}
and the beginning of the proof of Theorem \ref{MainTheorem}
in \S \ref{SectionProofMainTheorem}).

For the local $L$-functions,
we compare the $l$-adic and the $p$-adic case
(The following argument is inspired by \cite{TSaito}).
Fix a prime number $l \neq p$.
Since the weight-monodromy conjecture holds
for both $l$-adic and $p$-adic cohomology
(Theorem \ref{MainTheorem},
 Theorem \ref{MainTheoremp-adic}),
it is enough to prove the equality:
\begin{equation}
\label{EqualityTrace}
   \text{\rm Tr} \big( \Fr_q^m\,;\,H^{k}_{\text{\rm \'et}}(X_{\Gamma}
      \otimes_K \overline{K},\,\Q_l) \big)
    = \text{\rm Tr} \big( \varphi^m\,;\,D_{\text{\rm st}}(
        H^{k}_{\text{\rm \'et}}(X_{\Gamma}
        \otimes_K \overline{K},\,\Q_p)) \big)
\end{equation}
for all $k, m \geq 0$ (\cite{Rapoport}, Lemma 2.12).
Since 
$$ \dim_{\Q_l} H^{k}_{\text{\rm \'et}}(X_{\Gamma}
      \otimes_K \overline{K},\,\Q_l)
   = \dim_{\Q_p} H^{k}_{\text{\rm \'et}}(X_{\Gamma}
      \otimes_K \overline{K},\,\Q_p), $$
$H^{k}_{\text{\rm \'et}}(X_{\Gamma}
   \otimes_K \overline{K},\,\Q_p) = 0$
if $k$ is odd and $k \neq d$,
and $H^{k}_{\text{\rm \'et}} \big( X_{\Gamma}
      \otimes_K \overline{K},\,\Q_p(\frac{k}{2}) \big)$
is one dimensional and generated by algebraic cycles
if $k$ is even and $k \neq d$.
Therefore, (\ref{EqualityTrace}) trivially holds for $k \neq d$.
On the other hand,
the alternating sum of (\ref{EqualityTrace}) for all $k$
holds by Ochiai's results (\cite{Ochiai}, Theorem D).
Hence (\ref{EqualityTrace}) holds also for $k=d$.
\end{proof}

\begin{remark}
If $\X_{\Gamma}$ is strictly semistable,
we can directly prove Theorem \ref{LocalZetaFunctionP-adic}
without using Ochiai's results as follows.
It is easy to see that $E_1^{i,\,j}, d_1^{i,\,j}, N$
of the $l$-adic and the $p$-adic weight spectral sequences
have the same $\Q$-structure coming from the groups of
algebraic cycles on the special fiber.
Hence the characteristic polynomials of $\Fr_q, \varphi$
acting on each $E_2^{i,\,j}$ are equal.
Since the weight-monodromy conjecture holds
for both $l$-adic and $p$-adic cohomology
(Theorem \ref{MainTheorem}, Theorem \ref{MainTheoremp-adic}),
we see that the characteristic
polynomials of $\Fr_q, \varphi$ on
each graded quotient of the monodromy filtration on
$H^k_{\text{\rm \'et}}(X_{\Gamma} \otimes_K \overline{K},\,\Q_l)$,
$D_{\text{\rm st}}(H^{k}_{\text{\rm \'et}}(X_{\Gamma}
        \otimes_K \overline{K},\,\Q_p))$
are equal for all $k$.
Hence we have Theorem \ref{LocalZetaFunctionP-adic}.
\end{remark}

\begin{remark}
\label{RemarkP-adic}
The proof above is somewhat indirect
compared with the $l$-adic case.
The reason why we have to make such a detour
is that Schneider-Stuhler's results on de Rham cohomology
are not enough to compute the $p$-adic local zeta functions
because we do not have the action of $\varphi$
on de Rham cohomology in \cite{Schneider-Stuhler}.
However, by a recent work of Gro\ss{}e-Kl\"onne,
we have the action of $\varphi$ on the de Rham cohomology
of $\widehat{\Omega}^d_{K}$ (\cite{Grosse-Kloenne}).
By using his results,
it seems possible to compute the $p$-adic local zeta functions
directly without using $l$-adic cohomology.
\end{remark}

\section{Appendix : Application to the Tate conjecture for varieties
admitting $p$-adic uniformization}
\label{SectionApplicationTateConjecture}

Here we give an application of
the results in the previous section
to the Tate conjecture for varieties over number fields
admitting $p$-adic uniformization.
Such an application was pointed out to the author by Michael Rapoport
and the author includes it here with his permission.
The results in this appendix are essentially due to him.

Let $K$ be a number field,
$\overline{K}$ an algebraic closure of $K$,
and $X$ a proper smooth scheme over $K$ of dimension $d \geq 1$
which is not necessarily geometrically connected.
We say {\em $X$ admits $p$-adic uniformization}
if there is a finite extension $L$ of $K$ inside $\overline{K}$,
and a finite place $w$ of $L$ such that
each connected component of $X \otimes_{K} L_w$ is isomorphic to
$X_{\Gamma}$ for a cocompact torsion free discrete subgroup
$\Gamma \subset \PGL_{d+1}(L_w)$,
where $L_w$ is the completion of $L$ at $w$.

As an application of the results in the previous section,
we prove the Tate conjecture for $X$ as follows.

\begin{theorem}[cf. \cite{Tate}, Conjecture 1]
\label{ApplicationTateConjecture}
Let $X$ be a proper smooth scheme over a number field $K$
of dimension $d \geq 1$ admitting $p$-adic uniformization.
Fix a prime number $l$.
Then, for any finite extension $K'$ of $K$ inside $\overline{K}$
and for any $k,\ 0 \leq k \leq d$,
the $\Q_l$-vector space
$$ H^{2k}_{\text{\rm \'et}}\big(X \otimes_K \overline{K},
      \,\Q_l(k)\big)^{\Gal(\overline{K}/K')} $$
is generated by algebraic cycles on $X$
of codimension $k$ defined over $K'$.
\end{theorem}

\begin{remark}
In fact, as we see in the proof below,
the above $\Q_l$-vector spaces
are generated by power of the canonical bundles
of connected components of $X \otimes_K K'$.
\end{remark}

\begin{proof}
Since it is enough to prove the assertion after replacing $K'$
by a finite extension of it, we may assume $K'$ contains $L$.
Fix a finite place $w'$ of $K'$ over $w$,
and take a connected component $X'$ of $X \otimes_{K} K'$.
Then, $X' \otimes_{K'} K'_{w'}$ is 
isomorphic to $X_{\Gamma} \otimes_{L_w} K'_{w'}$
for some $\Gamma \subset \PGL_{d+1}(L_w)$.
We put
$V := H^{2k}_{\text{\rm \'et}}\big(X' \otimes_{K'} \overline{K},
    \,\Q_l(k)\big)$.

We shall prove that $V^{\Gal(\overline{K}/K')} $
is generated by a power of the canonical bundle $\L$ of $X'$.
If $k \neq d/2$, the assertion trivially holds
because ${\mathcal L}$ is ample by Mustafin, Kurihara
(\cite{Mustafin}, \cite{Kurihara}), and
$V$ is one dimensional by Schneider-Stuhler
(\cite{Schneider-Stuhler}, see also \S \ref{SectionSchneiderStuhlerConjecture}).

If $d$ is even and $k = d/2$, we shall prove that
$V^{\Gal(\overline{K}/K')}$ is one dimensional.
We know that this space has dimension $\geq 1$
because it contains the subspace generated by a power of $\L$.
To prove the opposite inequality, we consider the subspace
$V^{\Gal(\overline{K'_{w'}}/K'_{w'})}$ of $V$
consisting of elements fixed by the action of
$\Gal(\overline{K'_{w'}}/K'_{w'})$.
Here we fix an embedding
$\overline{K} \hookrightarrow \overline{K'_{w'}}$
over $K'$ and consider
$\Gal(\overline{K'_{w'}}/K'_{w'})$ as a subgroup of $\Gal(\overline{K}/K')$.
Since $V^{\Gal(\overline{K'_{w'}}/K'_{w'})}$ is one dimensional
by the following lemma and
$V^{\Gal(\overline{K}/K')} \subset V^{\Gal(\overline{K'_{w'}}/K'_{w'})}$,
we see that $V^{\Gal(\overline{K}/K')}$ is also one dimensional.
\end{proof}

\begin{lemma}
\label{LemmaTateConjecture}
As in \S \ref{SectionP-adic},
let $K$ be a finite extension of $\Q_p$,
$\Gamma \subset \PGL_{d+1}(K)$ a cocompact torsion free discrete subgroup,
and $K'$ a finite extension of $K$ inside $\overline{K}$.
Fix a prime number $l$, which is not necessarily different from $p$.
Assume that $d$ is even.
Then, the $\Q_l$-vector space
  $$ H^{d}_{\text{\rm \'et}} \big( X_{\Gamma} \otimes_{K} \overline{K},
      \,\Q_l{\textstyle \left( \frac{d}{2} \right)} \big)^{\Gal(\overline{K}/K')} $$
is one dimensional.
\end{lemma}

\begin{proof}
We put $V := H^{d}_{\text{\rm \'et}}(X_{\Gamma} \otimes_K \overline{K},\,\Q_l)$
for simplicity. Let $\F_{q'}$ be the residue field of $K'$.

Firstly, we consider the case $l \neq p$.
As in \S \ref{SubsectionLocalZetaFunctions}, we see that
$$ V^{I_{K'}} =
     \Q_l^{\oplus \mu(\Gamma)} \oplus \Q_l\!\left(\textstyle -\frac{d}{2} \right). $$
Recall that the monodromy filtration is invariant under
a finite extension of the base field.
Therefore,
$$ \big( V{\textstyle \left( \frac{d}{2} \right)} \big)^{\Gal(\overline{K}/K')}
  = \big( V{\textstyle \left( \frac{d}{2} \right)}^{I_{K'}} \big)^{\Fr_{q'}=1}
  = \big( {{\Q_l\!\left(\textstyle \frac{d}{2} \right)}}^{\oplus \mu(\Gamma)} \oplus
        \Q_l \big)^{\Fr_{q'}=1}
  = \Q_l $$
is one dimensional,
where \lq\lq $\Fr_{q'}=1$" denotes the subspace consisting of elements
fixed by the action of $\Fr_{q'}$.

The case $l=p$ is similar but we work with $D_{\text{\rm st}}(V)$
instead of $V$ itself.
Let
$m := \dim_{\Q_p} \big( V{\textstyle \left( \frac{d}{2} \right)} \big)^{\Gal(\overline{K}/K')}$.
We know that $m \geq 1$.
There is a $\Gal(\overline{K}/K')$-equivariant injection
${{\Q_p\!\left(\textstyle -\frac{d}{2} \right)}}^{\oplus m} \hookrightarrow V$.
By taking Fontaine's functor $D_{\text{\rm st}}$
for $\Gal(\overline{K}/K')$-representations,
we have an injection
$D_{\text{\rm st}} \big( {\Q_p\!\left(\textstyle -\frac{d}{2} \right)} \big)^{\oplus m}
    \hookrightarrow D_{\text{\rm st}}(V)$,
which is compatible with the actions of $N,\varphi$.
Since $N=0$ on $D_{\text{\rm st}} \big( {\Q_p\!\left(\textstyle -\frac{d}{2} \right)} \big)$,
the image must be contained in $D_{\text{\rm st}}(V)^{N=0}$.
Since $\varphi$ acts on $D_{\text{\rm st}} \big( {\Q_p\!\left(\textstyle -\frac{d}{2} \right)} \big)$
via multiplication by $q^d$,
$D_{\text{\rm st}}(V)^{N=0}$ has a subspace of dimension $m$
where $\varphi$ acts via multiplication by $q^d$.
On the other hand, as we see in \S \ref{SubsectionP-adicAnalogueLocalZetaFunctions},
the dimension of such a subspace is at most $1$.
Hence we have $m = 1$.

Hence the assertion of Lemma \ref{LemmaTateConjecture} is proved,
and the proof of Theorem \ref{ApplicationTateConjecture} is complete.
\end{proof}

\begin{remark}
It is known that certain Shimura varieties admit $p$-adic uniformization
(for example, see \cite{Rapoport-Zink2}, Theorem 6.50, Corollary 6.51).
Therefore, we can apply Theorem \ref{ApplicationTateConjecture}
to such Shimura varieties.
Moreover, since the Hasse-Weil zeta functions of them
can be written in terms of automorphic $L$-functions,
it seems possible to deduce another form of the Tate conjecture
for them concerning the orders of poles of the Hasse-Weil zeta functions
(see \cite{Tate}, Conjecture 2).
\end{remark}

\begin{remark}
In the function field case, Theorem \ref{ApplicationTateConjecture}
was essentially obtained by Laumon-Rapoport-Stuhler
by the same method as above (\cite{LaumonRapoportStuhler}, Proposition 16.2).
Note that, in the function field case,
the weight-monodromy conjecture was already known
to hold by Deligne (\cite{WeilII}).
\end{remark}

\end{document}